\documentstyle[12pt,amssymb]{article}
\begin{document}

\noindent{\bf HAMILTONIAN TYPE OPERATORS IN REPRESENTATIONS\\
OF THE QUANTUM ALGEBRA $U_q({\rm su}_{1,1})$}
\bigskip

{\sf N. M. Atakishiyev${}^1$ and A. U. Klimyk${}^{1,2}$}

{\it ${}^1$Instituto de Matem\'aticas, UNAM, CP 62210 Cuernavaca, 
Morelos, M\'exico}

{\it ${}^2$Institute for Theoretical Physics, Kiev 03143, Ukraine}

\bigskip

\begin{abstract}

We study some classes of symmetric operators for the discrete series
representations of the quantum algebra $U_q({\rm su}_{1,1})$, which
may serve as Hamiltonians of various physical systems. The problem of
diagonalization of these operators (eigenfunctions, spectra, overlap
coefficients, etc.) is solved by expressing their overlap coefficients
in terms of the known families of $q$-orthogonal polynomials. We consider
both bounded and unbounded operators. In the latter case they are not
selfadjoint and have deficiency indices (1,1), which means that they
may have infinitely many selfadjoint extensions. We find possible sets of
point spectrum (which depend on the representation space under
consideration) for one of such symmetric operators by using the 
orthogonality relations for $q$-Laguerre polynomials. 
In another case we are led to new orthogonality relations for
${}_3 \phi_1$-hypergeometric polynomials. Many new realizations for 
the discrete series representations are constructed, which follow 
from the diagonalization of the operators considered. In particular,
a new system of orthogonal functions on a discrete set is shown to
emerge.
\end{abstract}

\bigskip

\noindent{\sf I. INTRODUCTION}
\bigskip

The wealthy theory of representations of the Lie group $SU(1,1)\simeq
SL(2,{\Bbb R})$ and its Lie algebra (see, for example, Refs. 1, 2, 
Chapter 7, and 3) has been extensively employed in various branches
of physics and mathematics. Representations of the Lie algebra ${\rm su}(1,1)$
have been particularly useful in studying the isotropic harmonic oscillator,
non-relativistic Coulomb problem, relativistic Schr\"{o}dinger equation, Dirac
equation with the Coulomb interaction, and so on. The Hamiltonian $H$ in the
interacting boson model is represented as a linear combination of the 
operators,
corresponding to generating elements $J^{\rm cl}_+, J^{\rm cl}_-, 
J^{\rm cl}_0$ of the Lie algebra
${\rm su}(1,1)$. For this reason, the diagonalization of representation
operators, corresponding to such linear combinations, is an important
problem.

Diagonalization of representation operators is also important for solving
various other problems. In particular, matrices of a transition between bases,
diagonalizing two representation operators, describe (on the representation
level) automorphisms of the group $SU(1,1)$, if those two operators correspond
to one-parameter subgroups (see Ref. 4 and references therein).
The most recent papers on the diagonalization of representation operators for
$SU(1,1)$ and their use for elucidating properties of special functions are
Refs. 5--10.

After the appearance of quantum groups and quantum algebras${}^{11-13}$, 
most problems of the representation theory for Lie groups
and Lie algebras were transferred to the representation theory of quantized
groups and algebras. This development is also very important from the point
of view of possible applications both in mathematics and in physics. In
particular, the diagonalization of representation operators for simplest
quantum groups and algebras (especially, such as $U_q({\rm su}_2)$ and
$U_q({\rm su}_{1,1})$) is of great significance. These results are essential
for deeper understanding the theory of special functions and orthogonal
polynomials.

It was shown in Ref. 14 that for irreducible representations of 
$U_q({\rm su}_{1,1})$
a closure of a representation operator may be not a selfadjoint operator
(whereas its classical counterpart is selfadjoint). In particular,
the operator $A= q^{J_0/4}(J_++ J_-)q^{J_0/4} + bq^{J_0}$ in the positive
discrete series representation of $U_q({\rm su}_{1,1})$ was diagonalized
for $\vert b\vert < 2/(q^{-1/2}-q^{1/2})$. Later on, it turned out that
this operator is of the great importance for studying harmonic analysis
on the quantum group $SU_q(1,1)$ and for elucidating properties of
$q$-orthogonal polynomials. Namely, convolution identities for
Al-Salam--Chihara polynomials (which are very closely related to the
$q$-Meixner--Pollazcek polynomials) and for Askey--Wilson polynomials
were derived from the quantum algebra $U_q({\rm su}_{1,1})$ in Ref. 6. The
linearization coefficients for a two-parameter family of Askey--Wilson
polynomials were obtained. In Ref. 7, the Poisson kernel for Al-Salam--Chihara
polynomials is derived. By using the operator $A$, a bilinear generating
function for Askey--Wilson polynomials is obtained${}^8$. Diagonalization
of representation operators for $U_q({\rm su}_{1,1})$ is utilized to
interpret Askey--Wilson polynomials and ${}_8W_7$ series as matrix elements
of representation operators in mixed bases${}^{15}$. Spectral analysis of 
operators of irreducible $*$-representations of 
$U_q({\rm su}_{1,1})$ is employed in Ref. 16
to develop harmonic analysis on $SU_q(1,1)$ (see also Refs. 17 and 18).

Diagonalization of representation operators for the quantum algebras
$U_q({\rm su}_2)$ and $U_q({\rm su}_{1,1})$ finds wide applications in
physics. For example, some models in quantum optics, such as Raman and
Brillouin scattering, parametric conversion and the interaction of
two-level atoms with a single-mode radiation field (Dicke model), can
be described by interaction Hamiltonians, which are representation
operators for $U_q({\rm su}_2)$ or $U_q({\rm su}_{1,1})$ (see, for example,
Ref. 19 and references therein).

A great interest to spectral analysis of the operators for the positive
discrete series of $U_q({\rm su}_{1,1})$ appears in the analysis on
noncommutative (quantum) spaces. For example, the Laplace operator and
the squared radius (together with the third operator, which serves as
the operator $J_0$) generate the algebra $U_q({\rm su}_{1,1})$ 
acting on the
space of polynomials on the $n$-dimensional Manin space or on the quantum
complex vector space (see Refs. 20--22). These operators realize irreducible
representations of the algebra $U_q({\rm su}_{1,1})$, which belong to
the positive discrete series and form a $q$-analogue of the oscillator
representations of the Lie algebra ${\rm su}_{1,1}$. To construct Hamiltonians
of physical systems, existing in the Manin space or in the quantum complex
vector space (for example, Hamiltonians for harmonic oscillators in these
spaces), one thus needs to deal with operators of the positive discrete
series representations of $U_q({\rm su}_{1,1})$. Consequently, the
diagonalization of these Hamiltonians reduces to the problem of diagonalization
of representation operators for the quantum algebra $U_q({\rm su}_{1,1})$.

In the present paper we study the diagonalization (eigenfunctions, spectra,
transition coefficients, etc.) of some classes of operators for the discrete
series representations of the quantum algebra $U_q({\rm su}_{1,1})$, related
to $q$-orthogonal polynomials. We restrict ourselves by the discrete series
representations of $U_q({\rm su}_{1,1})$, because just these representation
operators are often used as Hamiltonians of physical models and these
representations are related to the $q$-oscillator algebra. The present
paper is a continuation of our study of representation operators, undertaken
in Ref. 23 (note that this paper has an essential overlap with Refs.
6, 7 and 15).

The paper is organized as follows. In section 2 we assemble some definitions
and formulas on discrete series representations of the algebras 
${\rm su}_{1,1}$
and $U_q({\rm su}_{1,1})$, which are necessary in the sequel. In section 3
we diagonalize the representation operator $I_1$ with a bounded continuous
spectrum, that is a $q$-extension of the ${\rm su}_{1,1}$-operator
$J^{\rm cl}_0 - J^{\rm cl}_1$. Then we briefly describe eigenfunctions and
spectra of a more general class of representation operators $I_1^{(\varphi)}$.
In section 4 we consider isometry between the standard Hilbert spaces 
${\cal H}_l$, on which the discrete series representations are realized, 
and the Hilbert space,
on which the operator $I_1$ is a multiplication operator. The latter Hilbert
space is the space of square integrable functions on the spectrum of the
operator $I_1$. We realize the representation $T^+_l$ on this Hilbert space.
In section 5 we study representation operator with bounded discrete spectrum.
We find this spectrum explicitly. This diagonalization is also related to a 
new realization of the representations of $U_q({\rm su}_{1,1})$, which is 
considered in section 6. In section 7 we consider the unbounded operator $I_3$.
This operator is symmetric, but not selfadjoint. Its deficiency indices are
$(1,1)$ and this means that $I_3$ has infinitely many selfadjoint extensions.
We find possible sets of point spectrum for this symmetric operator by
using the orthogonality relations for $q$-Laguerre polynomials. 
These sets depend on the space, on which the corresponding
representation of the algebra $U_q({\rm su}_{1,1})$ is realized.
In section 8 we consider the diagonalization of a one-parameter family of 
unbounded symmetric representation operators. These operators also have 
deficiency indices $(1,1)$ and  we find their eigenfunctions. 
In order to find selfadjoint extensions of the operators $I_3$ and $I_4$,
it is necessary to know extremal orthogonality measures for the 
orthogonal polynomials, related to these operators. To the best of our
knowledge, these measures are not explicitly known. Calculation of their 
orthogonality relations is a very complicated problem, which goes 
beyond the scope of the present paper.

Note that our results are well suited for studying properties of
$q$-orthogonal polynomials by different algebraic methods (including 
the methods of Refs. 6--8 and 15). This topic will be
discussed separately.

Throughout the sequel we always assume that $q$ is a fixed positive number
such that $q<1$. We extensively use the theory of $q$-special functions and
notations of the standard $q$-analysis (see, for example, Refs. 24 and 25). 
In particular, we assume that
$$
[a]_q:= \frac{q^{a/2}-q^{-a/2}}{q^{1/2}-q^{-1/2}}\, ,    \eqno (1.1)
$$
where $a$ can be a number or an operator.
\bigskip

\noindent{\sf II. DISCRETE SERIES REPRESENTATIONS OF ${\rm su}_{1,1}$
AND $U_q({\rm su}_{1,1})$}
\medskip

\noindent{\bf A. The classical case}
\medskip

The classical Lie algebra ${\rm su}_{1,1}$ is generated by the elements
$J^{\rm cl}_0, J^{\rm cl}_1, J^{\rm cl}_2$, satisfying the relations
$$
[J^{\rm cl}_0, J^{\rm cl}_1] = {\rm i} J^{\rm cl}_2\,,
\qquad [J^{\rm cl}_1, J^{\rm cl}_2] = -{\rm i}J^{\rm cl}_0\,,
\qquad [J^{\rm cl}_2, J^{\rm cl}_0] = {\rm i} J^{\rm cl}_1\,.
$$
In terms of the raising and lowering operators
$J^{\rm cl}_{\pm} = J^{\rm cl}_1 \pm {\rm i} J^{\rm cl}_2 $
these commutation relations can be written as
$$
[J^{\rm cl}_0, J^{\rm cl}_{\pm}] = \pm J^{\rm cl}_{\pm}\,,
\qquad  [J^{\rm cl}_-, J^{\rm cl}_+] = 2J^{\rm cl}_0\,. \eqno(2.1)
$$

The discrete series representations $T^+_l$ of ${\rm su}_{1,1}$ with
lowest weights are given by a positive number $l$ and they are realized
on the spaces ${\cal L}_l$ of polynomials in $x$. The basis of ${\cal L}_l$
consists of the monomials
$$
f_n^l (x) = c_n^{l}  x^n \, , \qquad n = 0, 1, 2, 3,\cdots \, , \eqno(2.2)
$$
where $c_n^{l} = \{(2 l)_n/ n!\}^{1/2}$. Assuming that this basis
consists of orthonormal elements, one defines a scalar product on
${\cal L}_l$. The closure of ${\cal L}_l$ leads to a Hilbert space,
on which the representation $T^+_l$ acts. The operators $A$, acting
on spaces of functions in $x$, will be denoted by $A(x)$.

One can consider an explicit realization of the representation operators
$J^{\rm cl}_i$, $i = 0,1,2$, in terms of the first-order differential
operators:
$$
J^{\rm cl}_0(x) = x {d\over dx} + l \,, \qquad J^{\rm cl}_1(x)
= {1\over 2}(1 + x^2) {d\over dx} + lx \,, \qquad
J^{\rm cl}_2(x) = {{\rm i}\over 2} (1 - x^2) {d\over dx} - {\rm i} lx \,. 
\eqno(2.3)
$$
Then the action of the generators $J^{\rm cl}_0, J^{\rm cl}_{\pm}$ in the
standard, or canonical, basis (2.2), consisting of the eigenfunctions of
the operator $J^{\rm cl}_0$, is given by
$$
J^{\rm cl}_+(x)\,f^{l}_n(x) = \sqrt{(2l + n)(n+1)}\, f^{l}_{n+1}(x)\,,\quad
J^{\rm cl}_-(x)\,f^{l}_n(x) = \sqrt{(2l+n-1)n}\,f^{l}_{n-1}(x)\,,
$$  $$
J^{\rm cl}_0(x)\,f_n^{l}(x) = (l+n)\,f_n^{l}(x)\, .
$$
In some cases it is of interest (see, for example, Ref. 4) to consider
a basis, associated with the eigenfunctions of the selfadjoint operator
$J^{\rm cl}_0(x) - J^{\rm cl}_1(x)$, namely,
$$
[J^{\rm cl}_0(x) - J^{\rm cl}_1(x)]\, {\eta}^l_{\lambda}(x) =
\lambda \, {\eta}^l_{\lambda}(x)\,, \eqno(2.4)
$$  $$
{\eta}^l_{\lambda}(x) = (1-x)^{-2l}\,\exp \left( \frac{2\lambda
x}{x-1}\right) \,.  \eqno(2.5)
$$
Here is an outline of a proof of (2.5); as we shall see in the sequel,
one can treat the quantum case quite similarly. Look for the
eigenfunctions ${\eta}^l_{\lambda}(x)$ in the form of an expansion
$$
{\eta}^l_{\lambda}(x) = \sum_{n=0}^{\infty}\,a_n^l(\lambda)\,f_n^l(x) =
\sum_{n=0}^{\infty}\,b_n^l(\lambda)\,x^n\,,   \eqno(2.6)
$$
where $b_n^l(\lambda) := a_n^l(\lambda)\,c_n^l$. Then
$$
[J^{\rm cl}_0(x) - J^{\rm cl}_1(x)]\,{\eta}^l_{\lambda}(x) =
\sum_{n=0}^{\infty}\,a_n^l (\lambda)\,
[J^{\rm cl}_0(x) - J^{\rm cl}_1(x)]\,f_n^l(x)
$$
$$ =
\sum_{n=0}^{\infty}\,a_n^l(\lambda)\,\left\{(n+l)f_n^l(x) -
\frac{1}{2} \sqrt {(2l+n)(n+1)}\, f^l_{n+1}(x) -
\frac{1}{2}\,\sqrt{(2l+n-1)n}\,f^l_{n-1}(x)\right\} \,. \eqno(2.7)
$$
This means that (2.4) leads to a three-term recurrence relation
$$
2(n+l-\lambda)\,b_n^l(\lambda) = (n+1)\,b^l_{n+1}(\lambda) +
(2l+n-1)\,b^l_{n-1}(\lambda)
\eqno(2.8)
$$
for the coefficients $b_n^l(\lambda)$. Since (2.8) represents the
recurrence relation for the Laguerre polynomials $L_n^{(2l-1)}(2\lambda)$,
one obtains that
$$
\eta^l_{\lambda}(x) = \sum_{n=0}^{\infty}\,b_n^l(\lambda)\,x^n =
\sum_{n=0}^{\infty}\,L_n^{(2l-1)}(2\lambda)\,x^n \,.  \eqno(2.9)
$$
Formula (2.5) now follows immediately from (2.9) upon using the
generating function
$$
\sum_{n=0}^{\infty}\,L_n^{(\alpha)}(x)\,t^n =
(1-t)^{-\alpha-1}\exp \left( \frac{xt}{t-1}\right)              \eqno(2.10)
$$
for the Laguerre polynomials
$$
L^{(2l-1)}_n(x ):= {(2l)_n\over n!}\, {}_1F_1 (-n, \, 2l;\ x) =
{(2l)_n\over n!}\,\sum_{k=0}^n {(-n)_k\over (2l)_k}{x^k\over k!} \,. 
\eqno(2.11)
$$

Observe that since the Laguerre polynomials (2.11) satisfy  the
orthogonality relation
$$
\int _0^\infty e^{- x}\, x^{2l-1}\,L^{(2l-1)}_m(x)\,L^{(2l-1)}_n(x)dx
= \frac{(n+2l-1)!}{n!}\,\delta_{mn}\,,                  \eqno(2.12)
$$
the spectrum of the operator $J^{\rm cl}_0(x) - J^{\rm cl}_1(x)$ is
simple and it covers the interval $[0,\infty )$.
\medskip

\noindent{\bf B. The quantum case}
\medskip

The quantum algebra $U_q({\rm su}_{1,1})$ is defined as the associative
algebra, generated by the elements $J_+$, $J_-$, and $J_0$, which satisfy
the commutation relations
$$
[J_0,J_{\pm}]=\pm J_{\pm},\ \ \ \ [J_-,J_+] =
{q^{J_0}-q^{-J_0}\over q^{1/2}-q^{-1/2}} \equiv [2J_0]_q \,, \eqno(2.13)
$$
and the conjugation relations
$$
J_0^* = J_0\,,\ \ \ \ \ J_+^* = J_-\, .         \eqno(2.14)
$$
(Observe that here we have replaced $J_-$ by $-J_-$ in the usual
definition of the algebra $U_q({\rm sl}_{2})$.) For convenience,
in what follows we denote the algebra $U_q({\rm su}_{1,1})$ by
${\rm su}_q(1,1)$.

The Casimir element of the algebra ${\rm su}_q(1,1)$ is given by
the formula
$$
C_q:= [J_0 -1/2]_q^2 - J_+\, J_- = [J_0 + 1/2]_q^2 - J_-\, J_+ \, .
$$

We are interested in the discrete series representations of ${\rm su}_q(1,1)$
with lowest weights. These irreducible representations will be denoted by
$T^+_l$, where $l$ is a lowest weight, which can take any positive number
(see, for example, Ref. 26). These representations are obtained by deforming
the corresponding representations of the Lie algebra ${\rm su}_{1,1}$.

As in the classical case, the representation $T^+_{l}$ can be realized
on the space ${\cal L}_{l}$ of all polynomials in $x$. We choose a basis
for this space, consisting of the monomials
$$
f^{l}_n\equiv f^{l}_n(x;q) := c^{l}_n(q)\, x^n, \ \ \  n = 0,1,2,\cdots ,
\eqno(2.15)
$$
where
$$ c^{l}_0(q) = 1, \qquad c^{l}_n(q) = \prod_{k=1}^n\,{[2l+k-1]_q^{1/2}
\over [k]_q^{1/2}} = q^{(1-2l)n/4}{(q^{2l};q)_n^{1/2}\over
(q;q)_n^{1/2}}\,, \ \ n = 1,2,3, \cdots,                   \eqno(2.16)
$$
and $(a;q)_n=(1-a)(1-aq)\ldots (1-aq^{n-1})$. The representation $T^+_{l}$
is then realized by the operators
$$
J_0(x) =x{d\over dx} + l,\qquad J_{\pm}(x) = x^{\pm 1}[ J_0(x) \pm l]_q \,. 
\eqno (2.17)
$$
As a result of this realization, we have
$$
J_+(x)\,f^{l}_n(x;q) =\sqrt{[2l+n]_q \,[n+1]_q }\, f^{l}_{n+1}(x;q),
$$ $$
J_-(x)\, f^{l}_n(x;q) =\sqrt{[2l+n-1]_q \,[n]_q }\,f^{l}_{n-1}(x;q), 
\eqno(2.18)
$$  $$
J_0(x)\, f^{l}_n(x;q) = (l + n)\,f^{l}_n(x;q) .
$$
Obviously, these operators satisfy the commutation relations (2.13).
The basis functions $f^{l}_n(x;q)$ are eigenfunctions of the operators
$J_0(x)$ and $C_q(x)$: $C_q(x)\, f^{l}_n(x;q)= [l - 1/2]^2_q\,f^{l}_n(x;q)$.

We know that the discrete series representations $T^+_l$ can be realized
on a Hilbert space, on which the conjugation relations (2.14) are satisfied.
In order to obtain such a Hilbert space, we assume that the monomials
$f^{l}_n(x;q)$, $n=0,1,2,\cdots$, constitute an orthonormal basis for this
Hilbert space. This introduces a scalar product $\langle \cdot ,\cdot \rangle$
into the space ${\cal L}_l$. Then we close this space with respect to this
scalar product and obtain the Hilbert space, which will be denoted by
${\cal H}_l$. The Hilbert space ${\cal H}_l$ consists of functions (series)
$$
f(x)=\sum _{n=0}^\infty b_nf^l _n(x;q)=\sum _{n=0}^\infty
b_nc^l_n(q) x^n= \sum _{n=0}^\infty a_nx^n,
$$
where $a_n=b_nc^l_n(q)$. Since $\langle f^l_m ,f^l_n\rangle
=\delta_{mn}$ by definition, for $f(x)=\sum _{n=0}^\infty a_nx^n$
and ${\tilde f}(x)=\sum _{n=0}^\infty {\tilde a}_n\,x^n$ we have
$\langle f , {\tilde f}\rangle = \sum _{n=0}^\infty a_n \,{\tilde
a}_n/|c^l_n(q)|^2$, that is, the Hilbert space ${\cal H}_l$
consists of analytical functions $f(x)=\sum _{n=0}^\infty a_n\,x^n$,
such that
$$
\Vert f\Vert ^2 \equiv \sum _{n=0}^\infty |a_n/c_n^l(q)|^2 <
\infty .
$$

It is directly checked that for an arbitrary function $f(x)\in
{\cal H}_l$ we have
$$
q^{c x {d\over dx}} \, f(x) = f (q^c x).
$$
Therefore, taking into account formulas (2.17), we conclude that
$$
q^{J_0(x)/2} \, f(x) = q^{{1\over 2} (x {d\over dx} + l)} \,
f(x) = q^{l/2} \, f(q^{1/2} x)\, ,                     \eqno(2.19)
$$ $$
J_+(x) \, f(x) = {x\over q^{1/2} - q^{-1/2}} \, \left [ q^l f(q^{1/2}x)
- q^{-l} f(q^{-1/2} x) \right ]\, ,                    \eqno (2.20)
$$ $$
J_-(x) \, f(x) = {1\over (q^{1/2} - q^{-1/2}) x} \, \left[
f(q^{1/2}x) - f(q^{-1/2} x) \right] \, .               \eqno (2.21)
$$
These relations will be used in the sequel for finding eigenfunctions
of representation operators.

We shall study spectra, eigenfunctions and overlap functions for operators
of the representations $T^+_l$, which correspond to elements of the quantum
algebra ${\rm su}_q(1,1)$ of the form
$$
A:= a\,q^{b\,J_0}[\alpha J_+ + {\alpha}^*
J_-]\,q^{b\,J_0} + f(q^{J_0}),\ \ \
a,b\in {\Bbb R}, \ \ \alpha\in {\Bbb C},\ \ |\alpha |=1, \eqno (2.22)
$$
where $f$ is some function.
These operators are representable in the basis (2.15) by a Jacobi matrix.
This circumstance allows one to apply the theory of $q$-orthogonal polynomials
for diagonalizing these operators.

A study of operators of the type (2.22) in 
representations of the quantum algebra $U_q({\rm su}_2)$ was started by
T. Koornwinder${}^{27}$ (see also Refs. 28 and 29), who  
applied them to investigation of the Askey--Wilson polynomials.

For the appropriate choice of the function $f$, the
operators (2.22) are symmetric, so one can close them and thus obtain
closed operators. If $A$ is bounded (for some particular choice of the
constants $a,b$ and $\alpha$), then it is a selfadjoint operator.
If $A$ is not bounded, then its closure may give symmetric operator,
which is not selfadjoint (see Ref. 14, note that its classical counterpart
is a selfadjoint operator). Deficiency indices of such closed operator
in this case are $(1,1)$ and it has infinitely many selfadjoint extensions.
Then the corresponding overlap coefficients (which are orthogonal polynomials)
may have many orthogonality relations. (We deal with operators of this type
in sections 7 and 8.) Only orthogonality relations, corresponding 
to extremal orthogonality measures, lead to selfadjoint extensions of
the operator $A$.
\bigskip

\noindent{\sf III. REPRESENTATION OPERATORS WITH BOUNDED CONTINUOUS SPECTRA}
\medskip

In this section we are interested in the operator
$$
I_1:= \frac {a}2\,q^{J_0} - b - q^{J_0/4}\,J_1\,q^{J_0/4} = 
\frac {a}2\, q^{J_0}-b-
\frac 12 \, \left[q^{1/4}J_+ + q^{-1/4}J_-\right]\,q^{J_0/2}     \eqno (3.1)
$$
of the discrete series representation $T^+_l$, where
$$
a=(q^{1/4}+q^{-1/4})\,b, \ \ \ b=(q^{1/2}-q^{-1/2})^{-1}\,.
$$
The representation operator $(q^{1/4}J_+ + q^{-1/4}J_-)q^{J_0/2}$ is
bounded (see [14]). Since $J_0$ has the eigenvalues $m=l,l+1,l+2,\cdots$,\,
the operator $q^{J_0}$ is also bounded (recall that $0<q<1$). Thus, the
operator $I_1$ is bounded. It is easy to check that $I_1$ is a selfadjoint
operator since
$$
I_1\,f^l_k = \beta_k(q)\,f^l_k- \alpha_k(q)\,f^l_{k+1} - \alpha_{k-1}
(q)\,f^l_{k-1}\,,
$$
where
$$
\alpha_k(q):=\frac 12 \left\{q^{l+k+1/2}\,[2l+k]_q\,[k+1]_q\right\}^{1/2},
\ \ \
\beta_k(q) =\frac{(q^{1/4}+q^{-1/4})q^{l+k}-2}{2(q^{1/2}-q^{-1/2})}.
$$
The constants $a$ and $b$ in (3.1) are chosen in such a way that in the limit
as $q\to 1$ the operator $I_1$ reduces to the ${\rm su}_{1,1}$-operator
$J^{\rm cl}_0 - J^{\rm cl}_1$ (see formula (2.4)).

Eigenfunctions of the operator $I_1(x)$,
$$
I_1(x)\,\xi_\lambda^l(x;q) = \lambda (q) \, \xi_\lambda ^l(x;q)\,,\eqno (3.2)
$$
and its spectrum can be found exactly in the same way as in the case of
the classical operator $J^{\rm cl}_0(x)-J^{\rm cl}_1(x)$ (observe that
these functions do not belong to the Hilbert space ${\cal H}_l$, if a
spectrum of the operator $I_1(x)$ is continuous). Namely, one can expand
these functions into Taylor series in $x$,
$$
\xi_\lambda^l(x;q) = \sum_{n=0}^{\infty}\,a^l_n(q)\,f^l_n(x;q)=
\sum_{n=0}^{\infty}\,b^l_n(q)\,x^n\,, \qquad b^l_n(q):=c^l_n(q)\,a^l_n(q)\,, 
\eqno(3.3)
$$
with coefficients of the expansion $b^l_n(q)$, and show that the $b^l_n(q)$
are expressed in terms of the continuous $q$-Laguerre polynomials. Then the
orthogonality relation for these polynomials determines a spectrum of the
operator $I_1(x)$.

However, we illustrate in this section a more "direct" way of deriving
eigenfunctions of the operator $I_1(x)$ by using  the explicit realization
(2.17) for the generators $J_0(x)$ and $J_{\pm}(x)$ 
(see Refs. 23, 28 and 29]).
Indeed, from formulas (2.19)--(2.21) it follows that
$$
I_1(x)\, f(x) = \frac {b}{2x}\, q^{(2l-1)/4}\,[(1-q^{(1-2l)/4}\,x)^2\,f(x) -
(1 - q^{l/2}\,x)(1 - q^{(l+1)/2}\,x)\,f(qx)]          \eqno (3.3)
$$
for an arbitrary function $f(x)$. It is thus natural to look for the
eigenfunctions $ \xi_{\lambda}^{l}(x;q)$ of the operator $I_1(x)$ in the form
$$
\xi^{l}_\lambda (x;q) = {(\alpha x; q)_\infty \, (\beta x;q)_\infty \over
(\gamma x; q)_\infty (\delta x; q)_\infty}\, ,   \eqno (3.4)
$$
where $(a;q)_\infty  = \prod _{r=0}^\infty (1-aq^r)$. Since
$(a;q)_\infty =(1-a)(aq ; q)_\infty$, we have
$$
\xi^{l}_\lambda (q x;q) =  {(1 - \gamma x) (1 - \delta x)\over (1 -
\alpha x)(1 - \beta x)}\, \xi^{l}_\lambda (x; q) \, . \eqno (3.5)
$$
Substituting (3.4) and (3.5) into (3.3), one gets
$$
I_1(x)\, \xi^{l}_\lambda (x;q) = \frac{b}{2x}\,q^{(2l-1)/4}\,
\left \{(1 - q^{(1-2l)/4}\,x)^2 \right.
$$ $$
\left.-\,(1 - q^{l/2}\,x)\,(1 - q^{(l+1)/2}\,x)\,{(1-\gamma x)(1- \delta x)
\over (1- \alpha x)(1- \beta x)}\right\}\, \xi^{l}_\lambda (x;q).   
\eqno (3.6)
$$
This equation can be written as
$$
I_1(x)\,\xi^{l}_\lambda (x;q) = \frac{b}2 \, q^{(2l-1)/4}\,\frac{A x^3 + 
B x^2 + Cx + D}
{(1 - \alpha x)(1 - \beta x)} \, \xi^{l}_{\lambda} (x;q) , \eqno (3.7)
$$
where the constant coefficients $A,B,C$, and $D$ are equal to
$$
A = q^{1/2}(q^{-l} \alpha \beta - q^{l}\gamma \delta),
$$  $$
B = q^{1/2}[q^{l} (\gamma + \delta)-q^{-l}(\alpha + \beta)] + (1+q^{1/2})
q^{l/2}\gamma \delta -2q^{(1-2l)/4}\alpha \beta ,
$$ $$
C =\alpha \beta -\gamma \delta+2q^{(1-2l)/4}(\alpha+ \beta)
- (1+q^{1/2})q^{l/2}(\gamma + \delta) - q^{1/2}\,(q^l - q^{-l}),
$$ $$
D =\gamma +\delta -\alpha -\beta +(1+q^{1/2})q^{l/2}-2q^{(1-2l)/4}.
$$
It is clear from (3.7) that the $\xi^{l}_\lambda (x;q)$ is an eigenfunction
of the operator $I_1(x)$ if the factor in front of $\xi^{l}_\lambda (x;q)$
on the right-hand side of (3.7) does not depend on $x$. It is the case if
$$ A = 0 , \qquad B = \alpha \beta D , \qquad C = - (\alpha + \beta) D\,.
\eqno (3.8)
$$
Then eigenvalues of the operator $I_1$ on the right-hand side of (3.7)
will be equal to $\lambda = q^{(2l-1)/4}D/2(q^{1/2} - q^{-1/2})$.
Requirements (3.8) are equivalent to the following three relations
between the parameters $\alpha , \beta , \gamma , \delta $:
$$
\alpha \beta = q^{2l} \gamma \delta ,\ \ \  (q^{1/2 - l} -
\alpha \beta ) (\alpha + \beta) = (q^{l +1/2} -  \alpha \beta)
(\gamma + \delta) - (1+q^{1/2})q^{-l/2}(q^l - q^{-l})\alpha \beta ,
$$ $$
(q^{l} - q^{-l})( q^{1/2} -  q^{-l} \alpha \beta ) = [\alpha +
\beta - (1+q^{1/2})q^{l/2}](\gamma + \delta - \alpha - \beta) .
$$
{}From these relations it follows that
$$ \alpha = q^{l/2} , \ \ \ \beta = q^{(l+1)/2},\ \ \
\gamma = q^{(1-2l)/4}\,e^{{\rm i}\theta } ,\ \ \
\delta = q^{(1-2l)/4}\,e^{-{\rm i}\theta } ,
$$
where $\theta$ is an arbitrary angle. Consequently, the eigenfunctions
of the operator $ I_1(x)$ are equal to
$$
\xi^l_{\lambda}(x;q) = \frac {(q^{l/2} x; q)_\infty \,(q^{(l+1)/2}
x; q)_{\infty}} {(q^{(1-2l)/4}  e^{{\rm i}\theta} x; q)_\infty
\,(q^{(1-2l)/4} e^{-{\rm i}\theta}x; q)_\infty}
= \frac {(q^{l/2}x;q^{1/2})_{\infty}}{(q^{(1-2l)/4} e^{{\rm i}
\theta}x;q)_{\infty}\, (q^{(1-2l)/4}e^{-{\rm i}\theta}x;q)_{\infty}}\, ,
 \eqno (3.9)
$$
and the corresponding eigenvalues are
$$
\lambda(q) =\frac{1-\nu }{q^{-1/2}-q^{1/2}}\, ,
$$
where $\nu = \cos \theta$ is a real parameter, see below.

The eigenfunctions (3.9) are in fact the generating functions for
the continuous $q$-Laguerre polynomials
$$
P_n^{(\alpha)}(y|q)=\frac{(q^{(2\alpha +3)/4}e^{-{\rm i}\theta};q)_n}
{(q;q)_n} \, q^{(2\alpha +1)n/4}e^{{\rm i}n\theta}
$$   $$
\times {}_2\phi _1(q^{-n},\ q^{(2\alpha +1)/4}e^{{\rm i}\theta};\
q^{-n-(2\alpha -1)/4}e^{{\rm i}\theta}\, ;\  q;q^{-(2\alpha -1)n/4}
e^{-{\rm i}\theta}),
$$
where $y=\cos \theta$ and ${}_2\phi _1$ is the basic hypergeometric
function, defined by formula (1.2.14) in Ref. 24. In order to make
this evident, one needs to represent (3.9) in the form
$$
\xi^l_{\lambda}(x;q) =\frac{(q^{2l-1/2}ax;q)_\infty (q^{2l}ax;q)_\infty}
{(q^{l-1/4}e^{{\rm i}\theta} ax;q)_\infty (q^{l-1/4}
e^{-{\rm i}\theta}ax;q)_\infty},
\eqno(3.10)
$$
where $a = q^{(1-3l)/2}$. Consequently, due to formula (3.19.11)
in Ref. 30, the desired connection is
$$
\xi^l_{\lambda}(x;q)= \sum _{n=0}^{\infty} q^{n(1-3l)/2}\, 
P_n^{(2l-1)}(\cos\theta|q)\,x^n
= \sum _{n=0}^{\infty} \frac{q^{n(1-3l)/2}}{c^l_n}\,
P_n^{(2l-1)}(\cos\theta|q)\,f_n^l(x) ,
$$
where $\cos \theta = 1 - (q^{-1/2}-q^{1/2})\lambda$.

Thus, we have proved that the eigenfunctions $\xi^l_\lambda(x;q)$
are connected with the basis elements $f^l_n(x;q)$ by the formula
$$
\xi^l_\lambda(x;q)=\sum _{n=0}^\infty p_n(\lambda)f^l_n(x;q), \eqno(3.11)
$$
where the overlap coefficients $p_n(\lambda)$ are given by
$$
p_n(\lambda)=
\frac{q^{(1/4-l)n}(q;q)_n^{1/2}}{(q^{2l};q)_n^{1/2}}
P^{(2l-1)}_n(\nu |q)=
\frac{q^{(1/4-l)n}(q;q)_n^{1/2}}{(q^{2l};q)_n^{1/2}}
P^{(2l-1)}_n(1 {-} (q^{-1/2} {-} q^{1/2})\lambda |q) \eqno (3.12)
$$
and, for convenience, $\lambda (q)$ is denoted by $\lambda$.

To find the spectrum of $I_1(x)$ (that is, a range of the parameter $\nu$),
we take into account the following. The selfadjoint operator $I_1(x)$ is
represented by a Jacobi matrix in the basis $f^l_n(x;q)$, $n=0,1,2,\cdots$.
As is evident from (3.11), the eigenfunctions $\xi^{l}_\lambda(x;q)$ are
expanded in the basis elements $f^l_n(x;q)$ with the coefficients (3.12).
According to the results of Chapter VII in Ref. 31, these polynomials
$p_n(\lambda)$ are orthogonal with respect to some measure $d\mu (\lambda)$
(this measure is unique, up to a multiplicative constant, since the operator
$I_1(x)$ is bounded). The set (a subset of ${\Bbb R}$), on which the 
polynomials
are orthogonal, coincides with the spectrum of the operator $I_1(x)$ and $d\mu
(\lambda)$ determines the spectral measure of this operator; the spectrum
of $I_1(x)$ is simple (see Chapter VII in Ref. 31).

We thus remind the reader that the orthogonality relation for the
continuous $q$-Laguerre polynomials $P_n^{(2l-1)}(y|q)$ has the form
$$
{1\over 2\pi} \int_{-1}^1P_m^{(2l-1)}(y|q)P_n^{(2l-1)}(y|q)
\frac{w(y)dy}{\sqrt{1-y^2}}= \frac{(q^{2l};q)_n \,q^{(2l-1/2)n}}
{(q;q)_\infty (q^{2l};q)_\infty (q;q)_n}\, \delta_{mn},
$$
where
$$
w(y)=\left\vert \frac{(e^{{\rm i}\theta};q^{1/2})_\infty (-e^{{\rm
i}\theta};q^{1/2})_\infty}{(q^{l-1/4}e^{{\rm i}
\theta};q^{1/2})_\infty} \right\vert ^2,\ \ \ y = \cos \theta
$$
(see formula (3.19.2) in Ref. 30). Therefore, the orthogonality
relation for the overlap coefficients (3.12) is
$$
\int _0^{2q^{1/2}/(1-q)} p_m(\lambda)p_n(\lambda){\hat
w}(\lambda)d\lambda =\delta_{mn}, \eqno (3.13)
$$
where
$$
{\hat w}(\lambda)=\frac{1}{2\pi} (q;q)_\infty (q^{2l};q)_\infty
\sqrt{ \frac{1-q}{\lambda q^{1/2}}} 
 \frac{w(1-q^{-1/2}(1-q)\lambda )}
{\sqrt{2-q^{-1/2}(1-q)\lambda }} . \eqno (3.14)
$$
Consequently, $\nu \in [-1,1]$ and the spectrum of the operator $I_1(x)$
(that is, a range of $\lambda\equiv \lambda (\nu;q)$) coincides with the
finite interval $[0,\,2q^{1/2}/(1-q)]$. The spectrum is continuous and
simple. The continuity of the spectrum means that the eigenfunctions
$\xi^{l}_\lambda(x;q)$ do not belong to the Hilbert space ${\cal H}_l$.
They belong to the space of functionals on ${\cal L}_l$,
which can be considered as a space of generalized functions on ${\cal L}_l$.
We have thus proved the following theorem.
\medskip

{\bf Theorem 1.} {\it The selfadjoint operator $I_1(x)$ has the continuous
and simple spectrum, which covers the finite interval $[0,\,2q^{1/2}/(1-q)]$.
The eigenfunctions $\xi^l_\lambda(x;q)$ are explicitly given by (3.9) and
they are related to the basis (2.15) by formula (3.11).}
\medskip

As we remarked at the beginning of this section, the operator $I_1(x)$
represents a $q$-extension of the ${\rm su}_{1,1}$-operator $J^{\rm cl}_0(x)
- J^{\rm cl}_1(x)$. In the limit as $q\to 1$ the finite interval
$[0,\,2q^{1/2}/(1-q)]$ of the eigenvalues of the operator $I_1(x)$ extends
to the infinite interval $[0,\infty)$. So if one puts $\nu = q^{\mu}$, then
$$
\lim_{q\to 1}\lambda (q^{\mu};q) = \mu \,. \eqno(3.15)
$$

Besides, it is known that the continuous $q$-Laguerre polynomials
$P_n^{(\alpha)}(y|q)$ have the following limit property (see [30],
formula (5.19.1))
$$
\lim_{q\to 1}P_n^{(\alpha)}(q^{\lambda}|q) = L_n^{(\alpha)}
(2\lambda)\,.\eqno(3.16)
$$
Thus, the coefficients of the series expansion (3.11) in $x$ of the
eigenfunctions $\xi^l_\lambda(x;q), \, \nu = q^{\mu}$, coincide with
the coefficients of the corresponding expansion of the
${\rm su}_{1,1}$-eigenfunctions $\eta^l_{\mu}(x)$\, (see (2.9))
in the limit as $q\to 1$.

There exists another, more complicated, family of selfadjoint
operators, closely related to $I_1$. They are defined as
$$
I_1^{(\varphi)} := \frac {a}{2}\,q^{J_0} - b  - q^{J_0/4}\,
\left[ \cos \varphi \, J_1
- \sin \varphi \,J_2\right]\,q^{J_0/4}
$$
$$ = \frac {a}{2}\, q^{J_0} - b - \frac 12\,\left[q^{1/4}\,
e^{{\rm i}\varphi}\,J_+
+ q^{-1/4}\,e^{-{\rm i}\varphi}\,J_-\right]\,q^{J_0/2}\,, \eqno(3.17)
$$
where $0\le \varphi < 2\pi$, $J_\pm =J_1\pm {\rm i}J_2$ and $a$, $b$
are such as in (3.1). These operators are bounded and
selfadjoint. Repeating the same reasoning, as for the operator
$I_1$, we arrive at the following theorem.
\medskip

{\bf Theorem 2.} {\it The eigenfunctions of the operator 
$I_1^{(\varphi )}(x)$ are
$$
\xi^l_\lambda(e^{{\rm i}\varphi}x;q) = \frac {(q^{l/2}\,
e^{{\rm i}\varphi} x; q)_\infty
\,(q^{(l+1)/2}\,e^{{\rm i}\varphi}x;q)_{\infty}} {(q^{(1-2l)/4}\,
e^{{\rm i}(\theta +\varphi)}x;q)_\infty \,(q^{(1-2l)/4}\,
e^{-{\rm i}(\theta -\varphi)}x;q)_\infty}\,,\qquad \nu = \cos\theta\,,
$$
where $\lambda=(1-\nu)/(q^{-1/2}-q^{1/2})$. Its spectrum is simple and
covers the interval $[0,\,2q^{1/2}/(1-q)]$, and the corresponding eigenvalues
$\lambda (\nu;q)$ are the same as for the operator $I_1(x)$}.
\medskip

The operators $I_1^{(\varphi)}$ are $q$-extensions of 
${\rm su}_{1,1}$-family of operators $J^{\rm cl}_0 - \cos\varphi\,
J^{\rm cl}_1 + \sin\varphi\,J^{\rm cl}_2$.
\bigskip

\noindent{\sf IV. REALIZATION OF REPRESENTATIONS $T^+_l$,
RELATED TO THE OPERATOR $I_1$}
\medskip

Let ${\cal L}^2_{0 ,a}:= {\cal L}^2([0 ,a],{\hat w}(\lambda)d\lambda )$
be the Hilbert space of functions $F(\lambda )$ on the interval $[0 ,a]$,
$a=2q^{1/2}/(1-q)$, with the scalar product
$$
\langle F_1,F_2\rangle_a =\int^{a}_0 F_1(\lambda)\overline{F_2(\lambda)}
{\hat w}(\lambda)d\lambda  ,                               \eqno (4.1)
$$
where ${\hat w}(\lambda)$ is determined by (3.14). Let us construct a
realization of the representation $T^+_l$ on the space ${\cal L}^2_{0 ,a}$.
Our reasoning in this section is close to the Favard theorem (see, for
example, Ref. 24).

Due to Theorem 1.6 of Chapter VII in Ref. 31, the space ${\cal D}$ of all
polynomials in $\lambda$ is everywhere dense in the Hilbert space
${\cal L}^2_{0,a}$. Due to formula (3.13), the polynomials (3.12) form
an orthonormal basis of ${\cal L}^2_{0,a}$.

Let ${\cal H}_l$ be the Hilbert space from section 2 and
$f(x)=\sum _{n=0}^\infty a_nf^l_n(x)$ be an expansion of $f\in {\cal H}_l$
with respect to the orthonormal basis (2.15). With every $f\in {\cal H}_l$
we associate a function $F(\lambda)$ on the spectrum of the operator $I_1$,
such that
$$
F(\lambda )=\langle f(x),\xi^{l}_{\lambda}(x;q)\rangle \equiv
\sum_{n=0}^\infty
a_n \,\frac{q^{(1/4-l)n}(q;q)^{1/2}_n}{(q^{2l};q)^{1/2}_n}
\,P^{(2l-1)}_n(1-(q^{-1/2}-q^{1/2})\lambda |q)              \eqno (4.2)
$$
(we have taken into account the expansion, given by formula (3.11)). Note
that the function $\xi^{l}_{\lambda}(x;q)$ (considered as a function of $x$)
does not belong to ${\cal H}_l$. This means that we consider $\langle f(x),
\xi^{l}_{\lambda}(x;q)\rangle$ as a formal expression, 
analogous to an integral
transform with a kernel, which does not belong to the corresponding Hilbert
space. By expanding $f(x)$ and $\xi^{l}_{\lambda}(x;q)$ with respect to the
basis $f^l_n(x)$ (see (3.11)), we formally obtain the right-hand side of 
(4.2). Clearly, the sum in (4.2) has a strict sense at least for 
functions $f(x)$ with
a finite number of nonzero coefficients $a_n$. For $f(x)\in {\cal H}_l$ this
sum converges in the topology of the Hibert space ${\cal L}^2_{0 ,a}$.
\medskip

{\bf Proposition 1.} {\it The mapping $\Phi : f(x)\to F(\lambda)$, given
by formula (4.2), establishes an invertible isometry between the Hilbert
spaces ${\cal H}_l$ and ${\cal L}^2_{0 ,a}$.}
\medskip

{\it Proof.} If $f(x)=\sum _n a_nf^l_n(x)\in {\cal H}_l$, then
due to the orthogonality relation (3.13) for the continuous $q$-Laguerre
polynomials from (4.2), we have $\langle F,F\rangle_a = \sum_n |a_n|^2 =
\langle f,f\rangle$. Invertibility of the mapping $\Phi$ is evident.
Proposition is proved.

It is easy to see that the isometry $\Phi$ maps basis elements $f^l_n$
of ${\cal H}_l$ to the basis elements $p_n(\lambda)$ of ${\cal L}^2_{0 ,a}$
from (3.12), respectively.

The results of the previous section allow us to realize the discrete
series representation $T^+_l$ of ${\rm su}_q(1,1)$ on the Hilbert space
${\cal L}^2_{0,a}$. Taking into account the self-adjointness of the operator
$I_1(x)$ and the fact that $I_1(x) \xi^{l}_{\lambda}(x;q)=
\lambda \xi^{l}_{\lambda}(x;q)$, for functions (4.2) we formally obtain
$$
I_1 F(\lambda)=\langle I_1(x)f(x),\xi^{l}_{\lambda}(x;q)\rangle=
\langle f(x),I_1(x) \xi^{l}_{\lambda}(x;q)\rangle = \lambda F(\lambda). 
\eqno (4.3)
$$
For this reason, we define an action of $I_1$ on the Hilbert 
space ${\cal L}^2_{0,a}$ by the formula $I_1\,f(\lambda ) = 
\lambda\,f(\lambda)$. Therefore, for the basis
elements (3.12) we get
$$
I_1\,p_k(\lambda )=\lambda\, p_k(\lambda). \eqno (4.4)
$$

Let us show that the operator $I_1$ acts on the basis $p_k(\lambda )$,
$k=0,1,2,\cdots $, by the formula
$$
I_1\, p_k(\lambda) = \frac{q^{l+k}(q^{1/4}+q^{-1/4})-2}
{2(q^{1/2}-q^{-1/2})}\, p_k(\lambda) - a_k \, p_{k+1}(\lambda )- 
a_{k-1}\, p_{k-1} (\lambda ), \eqno (4.5)
$$
where
$$
a_k=\frac 12\,({q^{k+l+1/2}\,[k+1]_q\,[k+2l]_q})^{1/2},
$$
that is, $I_1$ acts upon the basis functions $p_k(\lambda)$ of the space
${\cal L}^2_{0,a}$ by the same formulas as $I_1(x)$ acts upon the basis
elements $f^l_k$ of the space ${\cal H}_l$. In order to prove formula (4.5),
we replace (according to (4.4)) the left-hand side by $\lambda\, p_k(\lambda)$
and substitute the expression (3.12) for $p_{k-1}(\lambda)$, $p_{k}(\lambda)$,
and $p_{k+1}(\lambda)$. After simple transformations we obtain from (4.5)
the recurrence relation (3.19.3) of Ref. 30 for the continuous
$q$-Laguerre polynomials. This proves the formula (4.5).

As in the case of formula (4.3), for the action of the operator $q^{-J_0}$
we have
$$
q^{-J_0}\, p_k(\lambda )=\langle q^{-J_0} f^l_k(x;q), 
\xi^{l}_{\lambda}(x;q)\rangle
=q^{-l-k}\langle \, f^l_k(x;q), \xi^{l}_{\lambda}(x;q)\rangle =
q^{-l-k}\,p_k(\lambda),
$$
where $p_k(\lambda)$ are the basis elements (3.12), that is,
$$
q^{-J_0}\,p_k(\lambda)=q^{-l-k}\,p_k(\lambda).
$$
This means that $q^{-J_0}$ acts on the basis $p_k(\lambda )$, $k=0,1,2,
\cdots $, by the same formula as it acts on the basis (2.15) of the 
space ${\cal H}_l$.

It is easy to see that the operators $I_1$ and $q^{-J_0}$ determine
uniquely all other operators of the representation $T^+_l$. Thus, we 
have obtained a realization
of $T^+_l$ on the Hilbert space $ {\cal L}^2_{0,a}$.

Now we show that the operator $q^{-J_0}$ acts on ${\cal L}^2_{0,a}$ by
the formula
$$
q^{-J_0}\,F(\lambda) = q^{-l} \left( 1-{1 \over 4{\hat w}(\lambda;q^{2l})}\,
D_q {\hat w}(\lambda;q^{2l+1})\,D_q\right)\,F(\lambda),        \eqno (4.6)
$$
where ${\hat w}(\lambda;q^{2l}):= {\hat w}(\lambda)$ (${\hat w}(\lambda)$
is given by (3.14)) and $D_q\,f(\lambda ):= {f(\lambda)-f(q\lambda) 
\over \lambda
 - q\lambda }$. In order to prove formula (4.6), we take into account the
$q$-difference equation (3.19.5) from Ref. 30 for the continuous $q$-Laguerre
polynomials $P_n(y):=P^{(2l-1)}_n(y|q)$. From this $q$-difference equation
it follows that
$$
q^{-n}\,{\hat w}(y,q^{2l})\,P_n(y) = {\hat w}(y,q^{2l})\,P_n(y)-
{1\over 4}\,D_q \,{\hat w}(y,q^{2l+1})\,D_q\,P_n(y).
$$
Since $q^{-J_0}\, p_k(\lambda)=
q^{-l-k}\, p_k(\lambda )$, formula (4.6) does hold for the basis
elements $p_k(\lambda )$ and consequently for all $F\in {\cal L}^2_{0,a}$.
\medskip

{\bf Theorem 3.} {\it Let ${\cal L}^2_{0,a}$ be the Hilbert space,
introduced above. Then the representation $T^+_l$ can be determined
on it. The formulas (4.4) and (4.6) give the action of the operators
$I_1$ and $q^{-J_0}$ on ${\cal L}^2_{0,a}$. The operators $J_{\pm}$
and $J_0$ act on the basis $p_n(\lambda)$, $n=0,1,2,\cdots $, of this
space as $J_0\,p_n(\lambda)=(l+n)\,p_n(\lambda)$ and}
$$
J_+ \,p_n(\lambda) = \sqrt{[2l+n]_q\,[n+1]_q}\,p_{n+1}(\lambda),\ \ \
J_- \,p_n(\lambda) = \sqrt{[2l+n-1]_q\,[n]_q}\,p_{n-1}(\lambda).
$$
\medskip

The results of this section allow to prove the following assertion.
\medskip

{\bf Proposition 2.} {\it Let $p_n(\lambda)$ be the polynomials (3.12).
Then in the space ${\cal H}_l$ we have}
$$
p_n(I_1(x))\,f^l_0 = f^l_n.                      \eqno (4.7)
$$
\medskip

{\it Proof.} The isometry $\Phi : {\cal H}_l\to {\cal L}^2_{0,a}$ maps
$f^l_0\equiv 1$ to $p_0(\lambda)\equiv 1$. By formula (4.4) we have
$I_1^r\,p_0\equiv I_1^r\,1=\lambda ^k$. Therefore, $p_n(I_1)\,p_0 = 
p_n(\lambda)$.
Applying the mapping $\Phi ^{-1}$ to this identity, one obtains the desired
relation (4.7). Proposition is proved.
\bigskip

\noindent{\sf V. REPRESENTATION OPERATORS WITH BOUNDED DISCRETE SPECTRA}
\medskip

In this section we consider the operators
$$
I^{(\psi)}_2=q^{3J_0/4}(e^{{\rm i}\psi} J_+ + e^{-{\rm i}\psi}J_-)q^{3J_0/4} -
\left([J_0-l]_q\, q^{l/2}+[J_0+l]_q\,q^{-l/2}\right) q^{3J_0/2}
$$   $$
=(q^{3/4}e^{{\rm i}\psi} J_+ + q^{-3/4}e^{-{\rm i}\psi}J_-
-[J_0-l]_q\, q^{l/2}-[J_0+l]_q\,q^{-l/2})\,q^{3J_0/2},
$$
where $0<\psi \le 2\pi$ (note that these operators depend on the
index $l$ of the representation $T^+_l$). They act on the basis
elements (2.15) by the formula
$$
I^{(\psi)}_2\,f^l_k(x;q) = e^{{\rm i}\psi}\, q^{3(l+k)/2 + 3/4}
\sqrt{[k+1]_q\, [2l+k]_q}\,f^l_{k+1}(x;q)
$$  $$
+ e^{-{\rm i}\psi} \,q^{3(l+k)/2 - 3/4}\sqrt{[k]_q \,[2l+k-1]_q}\, 
f^l_{k-1}(x;q)
$$  $$
-\, q^{3(l+k)/2}\,([k]_q \,q^{(l-1)/2} + [2l+k]_q\, q^{-(l-1)/2})\,
f^l_{k}(x;q). \eqno (5.1)
$$
By using this action it is easy to check that the $I^{(\psi)}_2$ are
bounded selfadjoint operators for any value of $\psi\in (0,2\pi ]$. 
For a fixed value of $\psi$ we look for
eigenfunctions  of the operator $I^{(\psi)}_2$,
$$
I^{(\psi)}_2 \chi^l_\lambda (x;q) = \lambda \,\chi^l_\lambda (x;q),
$$
in the form 
$$
\chi^l_\lambda (x;q):=\sum _{k=0}^\infty P_{k}(\lambda)f^l_k(x;q)\,.
$$
As in section 2A, the equation
$$
I^{(\psi)}_2\,\chi^l_\lambda (x;q) = \sum _{k=0}^\infty P_{k}(\lambda)
\,I^{(\psi)}_2f^l_k(x;q)=\lambda \sum _{k=0}^\infty P_{k}(\lambda)\,f^l_k(x;q)
$$
leads to the following recurrence relation for the polynomials $P_k(\lambda)$,
which after simple transformations can be written as
$$
-e^{-{\rm i}\psi}\,q^{k+l}\,[(1-q^{k+1})(1-q^{2l+k})]^{1/2}\,P_{k+1}
(\lambda) - e^{{\rm i}\psi}\,q^{k+l-1}\,[(1-q^{k})(1-q^{2l+k-1})]^{1/2}\, 
P_{k-1} (\lambda)
$$ $$
+ \,(q^k-q^{2k+2l}+q^{2l+k-1}-q^{2k+2l-1})\,P_k(\lambda) =
(1-q^{-1})\, \lambda \, P_k(\lambda).                   \eqno (5.2)
$$
Upon making the substitution
$$
P_k(\lambda)= e^{{\rm i}k\psi}\,\left(\frac{(q^{2l};q)_k}
{(q;q)_k}\right)^{1/2}\, q^{-lk}\,P'_k(\lambda)
$$
in this recurrence relation, one derives the equation
$$
- q^{k}\,(1-q^{2l+k})\,P'_{k+1}(\lambda)- q^{k+l-1}\,(1-q^{k})\,
P'_{k-1}(\lambda)
$$ $$
+ \,(q^k-q^{2k+2l}+ q^{2l+k-1}-q^{2k+2l-1})\,P'_k(\lambda) =
(1-q^{-1}) \,\lambda \,P'_k(\lambda).
$$
This is the recurrence relation for the little $q$-Laguerre (Wall)
polynomials
$$
p_k(q^y;q^{2l-1}|q)={}_2 \phi_1 (q^{-k},0;\ q^{2l};\ q;q^{y+1})
$$ $$
=(q^{1-2l-k};q)_k\; {}_2 \phi_0 (q^{-k},q^{-y};\ -\, ;\ q;q^{y+-2l+1})
$$
with $q^y=(1-q^{-1})\lambda$. Thus, we have
$$
P'_k(\lambda)=p_k(q^y;q^{2l-1}|q),\ \ \ \ q^y=(1-q^{-1})\lambda,
$$
and, consequently,
$$
P_k(\lambda)=e^{{\rm i}k\psi}\left(\frac{(q^{2l};q)_k}{(q;q)_k}\right)^{1/2}
\,q^{-lk}\,p_k(q^y;q^{2l-1}|q). \eqno (5.3)
$$
This means that eigenfunctions of the operator $I^{(\psi)}_2$ are
of the form
$$
\chi ^l_\lambda (x;q) = c_l\,\sum _{k=0}^\infty q^{(1/4-3l)k}\,
e^{{\rm i}k\psi} \frac{(q^{2l};q)_k}{(q;q)_k}\,p_k(q^y;q^{2l-1}|q)\,
x^k,\ \ \ \ q^y=(1-q^{-1})\lambda.
$$

To sum up the right-hand side of this relation, one needs to know a
generating function
$$
F(x;\, t;\, a|q):= \sum_{n=0}^\infty \frac{(aq;q)_n}{(q;q)_n}\,
p_n(x;\, a|q)\,t^n  \eqno (5.4)
$$
for the little $q$-Laguerre polynomials
$$
p_n(x;\, a|q):={}_2 \phi_1 (q^{-n},0;\ aq;\ q;qx)=\frac{1}{(a^{-1}q^{-n};q)_n}
{}_2 \phi_0 (q^{-n},x^{-1};\ - ;\ q;x/a).               \eqno (5.5)
$$
To evaluate (5.4), we start with the second expression in (5.5) in terms
of the basic hypergeometric series ${}_2 \phi_0 $. Substituting it into
(5.4) and using the relation
$$
\frac{(q^{-n};q)_k}{(q;q)_k}=(-1)^kq^{-kn+k(k-1)/2}
\frac{(q;q)_n}{(q;q)_k(q;q)_{n-k}},
$$
one obtains that
$$
F(x;\, t;\, a|q):=\sum_{n=0}^\infty (-aqt)^nq^{n(n-1)/2}
\sum _{k=0}^n \frac{(x^{-1};q)_k}{(q;q)_k(q;q)_{n-k}}
(q^{-n}x/a)^k . \eqno (5.6)
$$
Interchanging the order of summations in (5.6) leads to the desired
expression
$$
F(x;\, t;\, a|q):=E_q(-aqt)\, {}_2 \phi_0 (x^{-1},0;\ - ;\ q;xt),    
\eqno (5.7)
$$
where $E_q(z)=(-z;q)_\infty$ is the $q$-exponential function of Jackson.

Similarly, if one substitutes into (5.4) the explicit form of the
little $q$-Laguerre polynomials in terms of ${}_2 \phi_1$ from
(5.5), this yields an expression
$$
F(x;\, t;\, a|q):=\frac{E_q(-aqt)}{E_q(-t)}\,{}_2 \phi_1 (0,0;\ q/t ;\ q;qx). 
\eqno (5.8)
$$

Using the explicit form of the generating function (5.7) for the little
$q$-Laguerre polynomials, we arrive at
$$
\chi ^l_\lambda (x;q)= (-e^{{\rm i}\psi}q^{(2l-1)/4}x;q)_\infty \,
{}_2\phi_0 (q^{-y},\, 0;\, - \,  ;\ q;\, e^{{\rm i}\psi}q^{y-(6l-1)/4}x), 
\eqno (5.9)
$$
where, as before, $q^y=(1-q^{-1})\lambda$. Another expression for
$\chi ^l_\lambda (x;q)$ can be written by using formula (5.8).

Due to the orthogonality relation
$$
(q^{2l};q)_\infty \sum_{k=0}^\infty \frac{q^{2lk}}{(q;q)_k}\,
p_m(q^k;q^{2l-1}|q) \,p_n(q^k;q^{2l-1}|q)= 
\frac{q^{2ln}(q;q)_n}{(q^{2l};q)_n}\, \delta_{mn}  \eqno (5.10)
$$
for little $q$-Laguerre polynomials (see formula (3.20.2) in 
Ref. 30), spectrum
of the operator $I^{(\psi)}_2$ coincides with the set of points
$q^n/(1-q^{-1})$, $n=0,1,2,\cdots$. This means that the eigenfunctions
$$
\chi ^l_{\lambda_n} (x;q)\equiv\Xi_n^l(x),\ \ \ \ n=0,1,2,\cdots , \ \ \
\lambda_n=\frac{q^n}{1-q^{-1}} ,                            \eqno (5.11)
$$
constitute a basis of the representation space. We thus proved the
following theorem.
\medskip

{\bf Theorem 4.} {\it The operator $I^{(\psi)}_2$ has a simple discrete
spectrum, which consists of the points $q^n/(1-q^{-1})$, $n=0,1,2,\cdots$.
The corresponding eigenfunctions $\Xi_n^l(x)$ constitute an orthogonal basis
in the space ${\cal H}_l$.}
\medskip

The basis (5.11) is orthogonal, but not orthonormal. The functions
$$
\hat\Xi ^l_n (x)=c_n\, \Xi ^l_n (x),\ \ \
c_n=q^{ln}\left( \frac{(q^{2l};q)_\infty}{(q;q)_n}\right) ^{1/2},
\ \ \ n=0,1,2,\cdots ,  \eqno (5.12)
$$
form an orthonormal basis of ${\cal H}_l$. This follows from the fact
that a matrix
$(a^l_{mn})$ with the entries $a^l_{mn}=c_nP_m(\lambda_n)$, which connects
the bases $\{ f^l_m\}$ and $\{ \hat\Xi ^l_n (x)\}$, is unitary (due to
the orthogonality relation for the little $q$-Laguerre polynomials).

The fact that the matrix $(a^l_{mn})$ is unitary and real means that
$$
\sum _{n=0}^\infty a^l_{mn}a^l_{m'n}=\delta _{mm'},\ \ \
\sum _{m=0}^\infty a^l_{mn}a^l_{mn'}=\delta _{nn'}.
$$
So, if one takes into account the explicit expression for $a^l_{mn}$, then
the first relation is actually the orthogonality relation for the little
$q$-Laguerre polynomials. The second relation corresponds to the 
orthogonality relation for the polynomials
$$
{\hat p}_n(q^{-m};\, q^{2l-1}| q):= (q^{-2l-m-1};q)_m\,
{}_2\phi_0(q^{-m},\, q^{-n};\, -\, :q,q^{-2l+1}q^n).
\eqno (5.13)
$$
This relation has the form
$$
\sum_{m=0}^\infty  q^{-lm} \frac{(q^{2l};q)_m}{(q;q)_m}
{\hat p}_n(q^{-m};\, q^{2l-1}| q)
{\hat p}_{n'}(q^{-m};\, q^{2l-1}| q)
=q^{-2ln} \frac{(q;q)_n}{(q^{2l};q)_n} \delta _{nn'}.
\eqno (5.14)
$$
Comparing (5.13) with the polynomials (3.15.1) in Ref. 30, we see 
that these polynomials are multiple to the Al-Salam--Carlitz polynomials
$V^{(2l-1)}_n(q^{-m};q)$ and (5.14) is the orthogonality
relation for them (compare with the relation (3.25.5) in Ref. 30). 
This means that the Al-Salam--Carlitz polynomials
$V^{(2l-1)}_n(q^{-m};q)$ are in fact dual with respect to the 
little $q$-Laguerre
polynomials. Thus, we see that our study of the operator
$I_2^{(0)}$ with the aid of the little $q$-Laguerre polynomials led us
to the orthogonality relation for the Al-Salam--Carlitz polynomials
$V^{(2l-1)}_n(q^{-m};q)$.
\bigskip

\noindent{\sf VI. A REALIZATION OF $T^+_l$, RELATED TO THE OPERATOR 
$I^{(\psi)}_2$}
\medskip

We know that the operator $I^{(\psi)}_2$ acts upon the basis
functions $\hat\Xi_n^l(x)$ as
$$
I^{(\psi)}_2 \hat\Xi_n^l(x)=\frac{q^n}{1-q^{-1}}\hat\Xi_n^l(x),\ \ \ \
n=0,1,2\cdots .  \eqno (6.1)
$$
We can also find how the operator $q^{-J_0}$ acts upon the basis (5.12).
{}From $q$-difference equation (3.20.5) 
in Ref. 30 for the little $q$-Laguerre
polynomials, one readily derives the following difference relation for the
polynomials $p_k(q^y)\equiv p_k(q^y;q^{2l-1}|q)$:
$$
q^{-k-l}p_k(q^y)=-q^{l-y-1}p_k(q^{y+1})+q^{-l}(1-q^{-y})p_k(q^{y-1})+
q^{-y}(q^{l-1}+q^{-l})p_k(q^y).
$$
As in the case of the representations of the algebra ${\rm su}_q(2)$
(see Ref. 29), we derive from it the action formula
$$
q^{-J_0}\Xi_n^l(x)=-q^{l-n-1}\Xi_{n+1}^l(x)+q^{-l}(1-q^{-n})
\Xi_{n-1}^l(x)+ q^{-n}(q^{l-1}+q^{-l})\Xi_n^l(x).
$$
In the orthonormal basis (5.12) this formula takes the form
$$
q^{-J_0}\hat\Xi_n^l(x)=q^{-n}\left[ (q^{l-1}{+}q^{-l})\hat\Xi_n^l(x) -
 q^{-1}(1{-}q^{n+1})^{1/2}\hat\Xi_{n+1}^l(x) - 
(1{-}q^{n})^{1/2} \hat\Xi_{n-1}^l(x)\right] .          
\eqno (6.2)
$$

The operators (6.1) and (6.2) completely determine the action of other
operators of the representation $T^+_l$ on the basis (5.12). However, the
corresponding formulas are not simple and we do not present them here.

Now we introduce the Hilbert space ${\frak l}^2_0$, which consists of
infinite sequences ${\bf a}=\{ a_k | k\in {\Bbb Z}_+ \}$ (note that
${\Bbb Z}_+=\{ 0,1,2,\cdots \}$), such that
$\sum_{k=0}^\infty q^{2lk}|a_k|^2/ (q;q)_k <\infty$.
The scalar product in this Hilbert space is naturally defined as
$$
\langle {\bf a},{\bf a}'\rangle_0=(q^{2l};q)_\infty \sum_{k=0}^\infty
\frac{q^{2lk}}{(q;q)_k} a_k\overline{a'_k} .
$$
Then by (5.10) the sequences of values of the polynomials from (5.3),
$$
P_n(\lambda)=\left( \frac{(q^{2l};q)_n)}{(q;q)_n}\right) ^{1/2}
q^{-ln} p_n((1-q^{-1})\lambda ;q^{2l-1}|q) ,
$$
on the set $\{ \lambda_k=\frac{q^k}{1-q^{-1}} |
k\in {\Bbb Z}_+ \}$ form an orthonormal basis of ${\frak l}^2_0$.
We denote these sequences by $\{ P_n(\lambda _k)\, |\,
k\in {\Bbb Z}_+ \}$.

Let ${\cal H}_l$ be the Hilbert space from section 2 and
$f(x)=\sum _{n=0}^\infty a_nf^l_n(x)$ be an expansion of $f\in {\cal H}_l$
with respect to the orthonormal basis (2.15). As in section 4, with every
function $f\in {\cal H}_l$ we associate the sequence $\{ F(\lambda_k)\, |\,
k\in {\Bbb Z}_+ \}$, $\lambda_k=\frac{q^k}{1-q^{-1}}$, such that
$$
F(\lambda _k)=\langle f(x),\eta^{l}_{\lambda_k}(x)\rangle ,
$$
where $\langle \cdot ,\cdot \rangle$ is the scalar product in
${\cal H}_l$. This defines a linear mapping $\Phi :f(x)\to 
\{ F(\lambda_k)\, |\,
k\in {\Bbb Z}_+ \}$ from ${\cal H}_l$ to the Hilbert space 
${\frak l}^2_0$. The following proposition is proved in the 
same manner as Proposition 1.
\medskip

{\bf Proposition 3.} {\it The mapping $\Phi : f(x)\to \{F(\lambda)\, |\,
k\in {\Bbb Z}_+ \}$ establishes an invertible isometry between the Hilbert
spaces ${\cal H}_l$ and ${\frak l}^2_0$.}
\medskip

It is easy to see that the isometry $\Phi$ maps basis elements $f^l_n$
of the space ${\cal H}_l$ to the basis elements $\{P_n(\lambda_k)\, |\,
k\in {\Bbb Z}_+ \}$ of the space ${\frak l}^2_0$.

As in section 4, we obtain a realization of the representation $T^+_l$
on the space ${\frak l}^2_0$, such that
$$
I_2^{(\psi)}\{ F(\lambda_k )\, |\, k\in {\Bbb Z}_+ \} =\{ \lambda_k
F(\lambda_k )\, |\, k\in {\Bbb Z}_+ \} ,\ \ \ \lambda_k = 
\frac{q^k}{1-q^{-1}} . \eqno (6.3)
$$
In particular, $I_2^{(0)}\{ P_n(\lambda_k )\, |\, k\in {\Bbb Z}_+ \}
=\{ \lambda_k P_n(\lambda_k)\, |\, k\in {\Bbb Z}_+ \} $. The operator
$I_2^{(0)}$ acts on $k$-th coordinate of the sequence $\{P_n(\lambda_k )\,
|\, k\in {\Bbb Z}_+ \} $ as $I_2^{(0)}P_n(\lambda_k ) = \lambda_k 
P_n(\lambda_k)$.
Taking into account formula (5.2), we deduce that
$$
I^{(0)}_2\{ P_n(\lambda_k )\, |\, k\in {\Bbb Z}_+ \}=
q^{3(l+n)/2}q^{3/4}([n+1]_q [2l+n]_q)^{1/2}\{ P_{n+1}(\lambda_{k}
)\, |\, k\in {\Bbb Z}_+ \}
$$  $$
+q^{3(l+n)/2}q^{-3/4}([n]_q [2l+n-1]_q)^{1/2}\{ P_{n-1}(\lambda_{k}
)\, |\, k\in {\Bbb Z}_+ \} 
$$  $$
-q^{3(l+n)/2}([n]_q q^{(l-1)/2}+[2l+n]_q q^{-(l-1)/2})\{
P_n(\lambda_k )\, |\, k\in {\Bbb Z}_+ \} , \eqno (6.4)
$$
that is, $I_2^{(0)}$ acts upon the basis elements $\{ P_n(\lambda)\,
|\, k\in {\Bbb Z}_+ \}$ by the same formulas as upon the basis
functions $f_n^l(x)$ of the space ${\cal H}_l$. We also have
$$
q^{-J_0}\{  P_n(\lambda_k )\, |\, k\in {\Bbb Z}_+ \} = q^{-l-n}
\{  P_n(\lambda_k )\, |\, k\in {\Bbb Z}_+ \} .     \eqno (6.5)
$$
The operators (6.4) and (6.5) determine uniquely all other operators
of the representation $T^+_l$ on ${\frak l}^2_0$. In particular,
we have
$$
J_+ \{P_n(\lambda_k )\, |\, k\in {\Bbb Z}_+ \} = \sqrt{[2l+n]_q\,[n+1]_q}
\,\{ P_{n+1}(\lambda_k )\, |\, k\in {\Bbb Z}_+ \} ,
$$  $$
J_-\{ P_n(\lambda_k )\, |\, k\in {\Bbb Z}_+ \} =
\sqrt{[2l+n-1]_q\,[n]_q}\,\{ P_{n-1}(\lambda_k )\, |\, k\in {\Bbb Z}_+ \}.
$$

The results of this section allow to prove (exactly in the same way as
Proposition 2) the following assertion.
\medskip

{\bf Proposition 4.} {\it Let $P_n(\lambda)$ be the polynomials
determined above. Then in the Hilbert space ${\cal H}_l$ we have}
$$
P_n(I_2^{(\psi)}(x))f^l_0=f^l_n.                  \eqno (6.6)
$$
\medskip

\noindent{\sf VII. UNBOUNDED REPRESENTATION OPERATORS}
\medskip

In this section we deal with the operator
$$
I_3:= - q^{-3J_0/4}(J_++J_-)q^{-3J_0/4} +
\frac{1+q}{2(1-q)}q^{-2J_0} - \frac{q^l+q^{1-l}}{2(1-q)}q^{-J_0}.
 \eqno (7.1)
$$
Since the operators $q^{-2J_0}$ and $q^{-J_0}$ are unbounded, the
$I_3$ is also an unbounded operator. We close the operator $I_3$ in
the space ${\cal H}_l$ and assume in what follows that $I_3$ is a
closed operator. It is easy to check that $I_3$ is a symmetric operator.
We shall see below that $I_3$ is not a selfadjoint operator and has
deficiency indices $(1,1)$, that is, it has infinite number of selfadjoint
extensions. 

In the limit as $q \to 1$ the operator $I_3$ reduces to the 
${\rm su}_{1,1}$-operator $J^{\rm cl}_0 - J^{\rm cl}_1 + 1/2$. 
In other words, (7.1)
represents another $q$-extension of essentially the same classical
operator $J^{\rm cl}_0 - J^{\rm cl}_1$ (see section 3).

To find eigenfunctions of the operator $I_3$,
$$
I_3(x)\,\zeta_\lambda ^l(x;q) = \lambda \,\zeta_\lambda^l(x;q),   \eqno (7.2)
$$
we evaluate first, by using (7.1) and (2.18), that
$$
I_3(x)\,f^l_n(x;q) =\frac{q^{-2(n+l)}}{2(1-q)}\{ \beta_n f^l_n(x;q)
- q^{-1/2}\,\alpha_n\,f^l_{n+1}(x;q) - q^{3/2}\alpha_{n-1} \,
f^l_{n-1}(x;q)\}. \eqno (7.3)
$$
where $\alpha_n=\sqrt{(1-q^{n+1})(1-q^{2l+n})}$ and $\beta_n=
1-q^{n+1}+q(1-q^{2l+n-1})$.
Substituting the expansion
$$
\zeta_\lambda^l(x;q):=\sum _{n=0}^\infty \, P_n(\lambda) \,f^l_n(x;q) =
\sum_{n=0}^{\infty}\,b_n^l\,x^n, \quad b^l_n:=P_n(\lambda)\,c_n^l(q),      
\eqno (7.4)
$$
into (7.2) and taking into account (7.3), one obtains the three-term
recurrence relation
$$
q^{-1/2}\sqrt{(1-q^{n+1})(1-q^{2l+n})}\,P_{n+1}^l(\lambda)
+ q^{3/2}\sqrt{(1-q^{n})(1-q^{2l+n-1})}\,P_{n-1}^l(\lambda)
$$    $$
= [1 - q^{n+1} + q(1 - q^{2l+n-1})(1+q)-
2(1-q)\,q^{2(l+n)}\lambda(q)]\,P_{n}^l(\lambda)        \eqno (7.5)
$$
for the coefficients $P^l_n(\lambda)$. Now multiplying both sides of (7.5) 
by $c^l_n(q)$ and using the relation $\sqrt{1-q^{2l+n}}\,c^l_n(q) =
q^{(2l-1)/4}\,\sqrt{1-q^{n+1}}\,c^l_{n+1}(q)$, one finally arrives at
the following recurrence relation for the coefficients $b^l_n$:
$$
q^{(2l-3)/4}\,(1-q^{n+1})\,b_{n+1}^l + q^{(7-2l)/4}\,(1-q^{2l+n-1})\, 
b_{n-1}^l
$$    $$
= \left[1-q^{n+1}+q\,(1-q^{2l+n-1})-2(1-q)\,
\lambda(q)\,q^{2(l+n)}\right] \,b_n^l.\eqno (7.6)
$$
Up to the multiplicative factor $q^{(3-2l)n/4}$, this is the recurrence
relation for $q$-Laguerre polynomials $L_n^{(\alpha)}(y;q)$, $\alpha=2l-1$,
defined as
$$
L_n^{(\alpha)}(y;q):=\frac{(q^{\alpha+1};q)_n}{(q;q)_n}\,
{}_1\phi _1 (q^{-n};\,q^{\alpha+1};\,q;\,-q^{n+\alpha+1}\,y) $$
$$ = \frac{(q^{\alpha+1};q)_n}{(q;q)_n}\,\sum_{k=0}^n \,
\frac{(q^{-n};q)_k\,y^k}{(q^{\alpha+1};q)_k\,(q;q)_k}\,
q^{k[n+\alpha+(k+1)/2]}\,.
$$
Consequently, we have
$$
b_n^l= q^{(3-2l)n/4}\,L_n^{(2l-1)}(\nu;q)\,,\ \ \
\lambda(q) = \lambda(\nu;q):=\frac{\nu}{2(1-q)}\, ,       \eqno (7.7)
$$
where $\nu$ is a real constant and the polynomials $P_n(\lambda)$ in
(7.4) are expressed in terms of the $q$-Laguerre polynomials as
$$
P_n(\lambda)=q^{n/2}\frac{(q;q)^{1/2}_n}{(q^{2l};q)^{1/2}_n}
L_n^{(2l-1)}(2(1-q)\lambda ;q).                         \eqno (7.8)
$$

This means that the eigenfunctions $\zeta_\lambda^l(x;q)$ of the
operator $I_3(x)$ have the form
$$
\zeta ^l_\lambda (x;q) = \sum _{n=0}^{\infty}\,b^l_n\,x^n
= \sum _{n=0}^{\infty}\,L_n^{(2l-1)}(\nu ;q)\,
\left(q^{(3-2l)/4}\,x\right)^n\,.                       \eqno (7.9)
$$
With the aid of the generating function
$$
\sum_{n=0}^{\infty}\,L_n^{\alpha}(x;q)\,t^n = \frac{1}{(t;q)_{\infty}}\,
{}_1\phi_1\,(-x;\,0;\,q;\,q^{\alpha+1}\,t)             \eqno (7.10)
$$
for the $q$-Laguerre polynomials (see formula (3.21.12) in [30]), one
can write these eigenfunctions as
$$
\zeta ^l_\lambda(x;q) = \frac{1}{(q^{(3-2l)/4}x;q)_{\infty}}
\, {}_1\phi _1(-\nu ;\,0;\,q; q^{3(1-2l)/4}\,x) ,        \eqno (7.11)
$$
where ${}_1\phi _1$ is the basic hypergeometric function.

Thus, the eigenfunctions $\zeta ^l_\lambda (x;q)$ can be expanded in
the basis (2.15) as
$$
\zeta ^l_\lambda (x;q) = \sum_{n=0}^\infty P_n(\lambda)\, f^l_n(x;q),
$$
where the polynomials $P_n(\lambda)$ are given by the formula (7.8).
As we can see from section 3.21 in Ref. 30, the polynomials $P_n(\lambda)$
have many different orthogonality relations. Taking into account the
results of Chapter VII in Ref. 31, we conclude from this fact that
{\it the closed symmetric operator $I_3(x)$ is not selfadjoint and
has deficiency indices} $(1,1)$. It has infinitely many selfadjoint
extensions. In order to find these extensions, it is necessary to know
extremal orthogonality measures for the $q$-Laguerre polynomials.
To the best of our knowledge, they are not known (see, for example,
Refs. 32 and 33).

The $q$-Laguerre polynomials have  
the discrete orthogonality relation
$$
b_c\sum _{k=-\infty}^\infty \frac{q^{2lk}}{(-cq^k;q)_\infty}
L_m^{(2l-1)}(cq^k;q)\,L_n^{(2l-1)}(cq^k;q)=
\frac{(q^{2l};q)_n}{(q;q)_nq^n} \delta _{mn} ,            \eqno (7.12)
$$
where $c$ is some positive number and
$$
b_c = c^{2l}\,q^{l(2l-1)}\ (q;q)_{2l-1}
$$
(see formula (3.21.3) in Ref. 30). Note that to each positive number
$c$ there corresponds an orthogonality relation.

We fix this positive number $c$ and introduce the Hilbert space
${\frak l}^2_c$, which consists of infinite sequences ${\bf a} =
\{ a_k | k\in {\Bbb Z}\}$, such that
$\sum _{k=-\infty }^\infty q^{2lk}|a_k|^2/ (-cq^k;q)_\infty
<\infty $.
The scalar product $\langle \cdot ,\cdot \rangle_c$ in ${\frak l}^2_c$
is given by
$$
\langle {\bf a} ,{\bf a}' \rangle_c =b_c \sum _{k=-\infty }^\infty
\frac{q^{2lk}}{(-cq^k;q)_\infty}\,a_k{\overline{a'_k}}.
$$
Then by (7.8) and (7.12) the sequences, consisting of values of the
polynomials $P_n(\lambda)$, $n=0,1,2,\cdots$, on the set
$\{\lambda_k = cq^k /(1-q)\,|\,k\in {\Bbb Z} \}$, form an
orthogonal system of elements in ${\frak l}^2_c$. We denote these
sequences by $\{P_n(\lambda _k)\, |\, k\in {\Bbb Z}\}$. However, as
we know from the results of Ref. 34, the polynomials
$\{P_n(\lambda _k)\, |\, k\in {\Bbb Z}\}$, $n=0,1,2,\cdots$, do not
constitute a basis in the space ${\frak l}^2_c$. We denote by
${\sf H}_c$ the closed subspace of ${\frak l}^2_c$, spanned by
the polynomials $\{P_n(\lambda _k)\, |\, k\in {\Bbb Z}\}$, $n=0,1,2,\cdots$.

Let ${\cal H}_l$ be the Hilbert space from section 2 and
$f(x)=\sum _{n=0}^\infty a_nf^l_n(x;q)$ be an expansion of $f\in{\cal H}_l$
with respect to the orthonormal basis (2.15). With every function
$f\in {\cal H}_l$ we associate the sequence $\{F(\lambda_k)\, |\, 
k\in {\Bbb Z}\}$, $\lambda_k = cq^k/2(1-q)$, such that
$$
F(\lambda _k)=\langle f(x),\zeta^{l}_{\lambda_k}(x;q)\rangle =
\sum_{n=0}^\infty a_n  P_n(\lambda_k ),                  \eqno (7.13)
$$
This yields the linear mapping $\Phi :f(x)\to \{ F(\lambda_k)\, |\,
k\in {\Bbb Z}\}$ from ${\cal H}_l$ to ${\sf H}_c$. The following
proposition is proved in the same way as Proposition 1.
\medskip

{\bf Proposition 5.} {\it The mapping $\Phi : f(x) \to \{F(\lambda_k)\,
|\, k\in {\Bbb Z}\}$ establishes an invertible isometry from ${\cal H}_l$
to the Hilbert spaces ${\sf H}_c$.}
\medskip

It is easy to see that the isometry $\Phi$ maps basis functions $f^l_n(x;q)$
of ${\cal H}_l$ to the elements $\{P_n(\lambda_k)\,|\, k\in {\Bbb Z}\}$
of the space ${\sf H}_c$.

As in section 4, we can check that the isometry $\Phi$ transforms the
operator $I_3(x)$ to the operator $\hat I_3$, which acts on the linear
span of the polynomials $P_n(\lambda)$, $n=0,1,2,\cdots$, as the
multiplication operator, $\hat I_3\,P_n(\lambda ) = \lambda \,P_n(\lambda)$,
if $P_n$ are considered on the set $\lambda_k=cq^k/2(1-q)$, $k\in {\Bbb Z}$.
We extend the operator $\hat I_3$ to the multiplication operator $I_3^{(c)}$
on the space ${\sf H}_c$:
$$
I_3^{(c)}\{ F(\lambda_k )\, |\, k\in {\Bbb Z}\} =\{ \lambda_k
F(\lambda_k )\, |\, k\in {\Bbb Z}\} ,\ \ \ \lambda_k = 
\frac{cq^k}{2(1-q)} . \eqno(7.14)
$$

As in the case of formula (4.5), with the aid of the recurrence relation for
the $q$-Laguerre polynomials one obtains that $I_3^{(c)}$ acts upon
the basis elements $\{ P_n(\lambda_k)\, |\, k\in {\Bbb Z}\}$ of
${\sf H}_c$ by the same formula as $I_3(x)$ acts upon the basis
functions $f_n^l(x;q)$ of the space ${\cal H}_l$, that is, by the
formula (7.3). To the element $q^{J_0}$ there corresponds the operator
$$
q^{J^{(c)}_0}\{ P_{n}(\lambda_k )\, |\, k\in {\Bbb Z}\} =
q^{l+n}\{P_{n}(\lambda_k )\, |\, k\in {\Bbb Z}\}
$$
on the space ${\sf H}_c$. The operators $I_3^{(c)}$ and $q^{J^{(c)}_0}$
determine uniquely all other operators of the representation $T^+_l$
on the space ${\sf H}_c$. Thus, we have proved the following theorem.
\medskip

{\bf Theorem 5.} {\it Let ${\sf H}_c$ be the Hilbert space, defined
above. Then the representation $T^+_l$ is defined on ${\sf H}_c$.
In particular, the operators $q^{J_0}$ and $J_\pm$ act on the
basis $\{ P_{n}(\lambda_k )\, |\, k\in {\Bbb Z}\}$ of this space as
$q^{J_0}\{ P_{n}(\lambda_k )\, |\, k\in {\Bbb Z}\} = q^{l+n}
\{P_{n}(\Lambda_k ) \, |\, k\in {\Bbb Z}\}$ and}
$$
J_+\{ P_n(\lambda_k )\, |\, k\in {\Bbb Z}_+ \} = \sqrt{[2l+n]_q\,[n+1]_q}
\{ P_{n+1}(\Lambda_k )\, |\, k\in {\Bbb Z}_+ \} ,
$$  $$
J_-\{ P_n(\lambda_k )\, |\, k\in {\Bbb Z}_+ \} = \sqrt{[2l+n-1]_q\,[n]_q}
\{ P_{n-1}(\lambda_k )\, |\, k\in {\Bbb Z}_+ \}.
$$
\medskip

The operator of multiplication by independent variable on the 
space ${\frak l}^2_c$ is selfadjoint (see section 48, 
Chapter 4, in Ref. 35) and its restriction to ${\sf H}_c$ 
gives the operator $I^{(c)}_3$. However,
the self-adjointness of the former operator does not mean that the
latter operator $I^{(c)}_3$ is selfadjoint. However, we know from
the orthogonality relation for the $q$-Laguerre
polynomials that the set of points
$$
cq^k/2(1-q) , \ \ \ \ k=0,\pm 1,\pm 2,\cdots ,
$$
constitute the point spectrum of $I^{(c)}_3$.

The functions
$$
\zeta^l_{\lambda_k}(x;q)\equiv \Omega ^{l,c}_k(x),\ \ \
\lambda _k=cq^k/2(1-q),\ \ \ k=0,\pm 1,\pm 2,\cdots ,
$$
defined by (7.11), constitute a basis of the Hilbert space ${\cal H}_l$.
This basis is orthogonal, but not orthonormal. The functions
$$
\hat\Omega ^{l,c}_k(x)=d_k \Omega ^{l,c}_k(x), \ \ \
k=0,\pm 1,\pm 2,\cdots ,                            \eqno (7.15)
$$
with $d_k=b^{1/2}_cq^{lk}/(-cq^k;q)^{1/2}_\infty$ form an
orthonormal basis in ${\cal H}_l$. This can be shown in the same way
as in the case of the basis (5.12).

We know that the operator $I_3(x)$ acts upon the basis functions
$\hat\Omega ^{l,c}_n(x)$ as
$$
I_3(x)\hat\Omega ^{l,c}_n(x)=\frac{cq^n}{2(1-q)}
\hat\Omega ^{l,c}_k(x).                             \eqno (7.16)
$$
We can find how the operator $q^{J_0}$ acts upon this basis.
{}From $q$-difference equation (3.21.6) in Ref. 30 for the $q$-Laguerre
polynomials one finds the following difference equation for the
polynomials $p_n(y)\equiv L_n^{2l-1}(y;q)$:
$$
q^n yp_n(y)=(1+y)p_n(qy)-(q^{-2l+1}+1)p_n(y)+q^{-2l+1}p_n(q^{-1}y).
$$
As in the case of the representations of the algebra ${\rm su}_q(2)$
(see Ref. 29), this leads to the action formula
$$
q^{J_0} \Omega ^{l,c}_k(x)=q^l(c^{-1}q^{-k}+1)\Omega ^{l,c}_{k+1}(x)
+c^{-1}q^{1-l-k}\Omega ^{l,c}_{k-1}(x)-
c^{-1}q^{l-k}(1+q^{-2l+1})\Omega ^{l,c}_{k}(x).
$$
In the orthonormal basis $\{\hat\Omega ^{l,c}_n(x) \}$ this formula
takes the form
$$
q^{J_0} \hat\Omega ^{l,c}_k(x)=\frac{1}{cq^k}
 \left[ (1+cq^{k})^{1/2}\hat\Omega ^{l,c}_{k+1}(x)
+q(1+cq^{k-1})^{1/2} \hat\Omega ^{l,c}_{k-1}(x)
-(q^l+q^{1-l})\hat\Omega ^{l,c}_{k}(x)\right] . \eqno (7.17)
$$
The formulas (7.16) and (7.17) completely determine other operators
of the representations $T^+_l$ in the basis $\{\hat\Omega ^{l,c}_n(x) \}$.

There exists another, more complicated, family of symmetric operators,
closely related to $I_3$. They are of the form
$$
I^{(\psi)}_3:=  - q^{-3J_0/4}( e^{{\rm i}\psi} J_+ + e^{-{\rm i}\psi}
J_-)q^{-3J_0/4} + \frac{1+q}{2(1-q)}\, q^{-2J_0} - \frac{q^l + 
q^{1-l}}{2(1-q)}\, q^{-J_0}
$$
and can be diagonalized in exactly the same way as above.

At the end of this section we mention that the results obtained above
lead to a new system of orthogonal functions on a discrete set. They can be
obtained from the entries of the matrix $(a_{kn})$, which connects the
orthonormal bases $f^l_n(x)$, $n=0,1,2,\cdots$, and 
$\hat\Omega ^{l,c}_k(x)$, $k=0,\pm 1,\pm 2,\cdots$. These entries
are of the form
$$
a_{kn}=d_kq^{n/2}\frac{(q;q)_n^{1/2}}{(q^{2l};q)_n^{1/2}}
L_n^{(2l-1)}(cq^k;q),
$$
where $d_k$ are defined in (7.15). Orthonormality of the both bases
means that
$$
\sum _{k=-\infty}^\infty a_{kn}a_{kn'}=\delta _{nn'},\ \ \
\sum _{n=0}^\infty a_{kn}a_{k'n}=\delta _{kk'} .
$$
The first relation in fact coincides with the orthogonality
relation (7.12) for the $q$-Laguerre polynomials. The second
relation corresponds to the orthogonality relations for the set of
functions
$$
F_k(q^{-n};\, c, 2l-1|q):=(q;q)_n\, {}_2 \phi_1(q^{-n},-cq^k;\,
0;\; q,q^{n+2l}),\ \ \ k=0,\pm 1,\pm 2,\cdots.
\eqno (7.18)
$$
This orthogonality relation is of the form
$$
\sum_{n=0}^\infty \frac{q^n(q;q)_n}{(q^{2l};q)_n}
F_k(q^{-n};\, c, 2l-1|q)F_{k'}(q^{-n};\, c, 2l-1|q)
=\frac{(-cq^k;q)_\infty}{ q^{2lk}b_c}\delta _{kk'}.
\eqno (7.19)
$$
Thus, in the study of the operator $I_3$ with the aid of the $q$-Laguerre
polynomials a new system of functions, orthogonal on a
discrete set, emerges.
\bigskip

\noindent {\sf VIII. ANOTHER EXAMPLE OF THE UNBOUNDED CASE}
\medskip

In this section we wish to find eigenfunctions of the symmetric
operator
$$
{I}^{(\psi)}_4=q^{-J_0/4}(e^{{\rm i}\psi}J_+ + e^{-{\rm i}\psi}J_-)q^{-J_0/4}
= (e^{{\rm i}\psi}q^{-1/4}J_++e^{-{\rm i}\psi}q^{1/4}J_-)q^{-J_0/2}.
$$
This operator acts upon the basis elements $f^l_k$ as
$$
{I}^{(\psi)}_4f^l_k=q^{-l-n-1/2}\left( e^{{\rm i}\psi}\sqrt{(1{-}q^{k+1})
(1{-}q^{2l+k})} f^l_{n+1}+q e^{-{\rm i}\psi}\sqrt{(1{-}q^{k})
(1{-}q^{2l+k-1})} f^l_{n-1}\right) .
$$
Since the coefficients here tend to $\infty$ when $n\to \infty$, the
operators ${I}^{(\psi)}_4$ are unbounded. Let us show that the closure
of $I^{(\psi)}_4$ is not a selfadjoint operator and has selfadjoint
extensions (actually, there are many such extensions). Changing the
basis $\{ f^l_k\}$ by the basis $\{ {\tilde f}^l_k\}$, where
${\tilde f}^l_k=e^{{\rm i}k\psi}f^l_k$, we obtain the matrix
form of the operator ${ I}^{(\psi)}_4$ in the new basis:
$$
{I}^{(\psi)}_4 {\tilde f}^l_k=a_k {\tilde f}^l_{k+1}+ a_{k-1}
{\tilde f}^l_{k-1}, \ \ \ \
a_k=q^{-l-k-1/2}\sqrt{(1-q^{k+1})(1-q^{2l+k})}.
$$
According to Theorem 1.5 in Chapter VII of Ref. 31, the closed
operator ${I}^{(\psi)}_4$ is not selfadjoint and has deficiency
indices (1,1) if $a_{k-1}a_{k+1}\le a^2_k$ and $\sum _{k=0}^\infty a_k^{-1}
<\infty$. Let us show that these conditions are fulfilled for the
operator ${I}^{(\psi)}_4$.

Since $q+q^{-1}\ge 2$ for positive values of $q$ (and $q+q^{-1}=2$
only if $q=1$), then
$$
(1-q^{k+1})(1-q^{k-1})\le (1-q^k)^2.
$$
This leads to $a_{k-1}a_{k+1}\le a^2_k$. Since
$a_k/a_{k+1} \to q<1$ when $k\to \infty$
and $0<q<1$, then $\sum _{k=0}^\infty a_k^{-1}<\infty$. Thus,
${ I}^{(\psi)}_4$ is not selfadjoint operator and has deficiency
indices (1,1).

For eigenfunctions $\eta _\lambda (x;q)$ of the operator ${I}^{(\psi)}_4(x)$
we have
$$
{I}^{(\psi)}_4 (x)\eta _\lambda (x;q)=\lambda \eta _\lambda (x;q)
$$
and these functions can be represented as
$$
\eta _\lambda (x;q)=\sum _{k=0}^\infty P_k(\lambda)f^l_k(x;q). \eqno (8.1)
$$
Then for polynomials $P_k(\lambda)$ we obtain the recurrence
relation, which is equivalent to the following one:
$$
e^{-{\rm i}\psi} \sqrt{(1-q^{k+1})(1-q^{2l+k})} P_{k+1}(\lambda)+
qe^{{\rm i}\psi} \sqrt{(1-q^{k})(1-q^{2l+k-1})} P_{k-1}(\lambda)
$$   $$
=q^{l+k-1/2}(q^{-1/2}-q^{1/2})\lambda P_k(\lambda). \eqno (8.2)
$$
Upon making the substitution
$$
P_k(\lambda)=e^{{\rm i}k\psi}\frac{(q;q)_k^{1/2}} {(q^{2l};q)_k^{1/2}}\,
 P'_k(\lambda)\,,                                     \eqno (8.3)
$$
one obtains
$$
(1-q^{k+1})P'_{k+1}(\lambda)+ q(1-q^{2l+k-1}) P'_{k-1}(\lambda)=
q^kd\lambda P'_k(\lambda) ,                            \eqno (8.4)
$$
where $d=(1-q)\,q^{l-1}$.

Now we substitute the expression (8.3) for $P_k(\lambda)$ and the
expression (2.15) for $f^l_k(x;q)$ into (8.1). After simple
transformation this yields
$$
\eta ^l_\lambda (x;q)=\sum _{k=0}^\infty P'_k(\lambda)
(q^{(1-2l)/4}e^{{\rm i}\psi}x)^k .                     \eqno (8.5)
$$

Multiply both sides of (8.4) by $y^{k+1}$, $y:=q^{(1-2l)/4}e^{{\rm i}\psi}x$,
and sum up over $k$:
$$
\eta ^l_\lambda (x;q)-\eta ^l_\lambda (qx;q)+y^2q\eta ^l_\lambda (x;q)-
y^2q^{2l+1}\eta ^l_\lambda (qx;q)=yd\lambda \eta ^l_\lambda (qx;q).
$$
This gives
$$
\eta ^l_\lambda (x;q)=\frac{1+yd\lambda +y^2q^{2l+1}}{1+qy^2} \eta
^l_\lambda (qx;q).                                      \eqno (8.6)
$$
Then we set $d\lambda =-2q^{l-1/2}\cosh \theta$ and iterate (8.6):
$$
\eta ^l_\lambda (x;q)=\frac{(yq^{l+1/2}e^\theta;q)_n (yq^{l+1/2}
e^{-\theta};q)_n}
{({\rm i}yq^{1/2};q)_n (-{\rm i}yq^{1/2};q)_n} \eta ^l_\lambda (q^nx;q).
$$
Passing to the limit $n\to \infty$ and taking into account that
$q^n\to 0$ when $n\to \infty$, we arrive at the following explicit 
expression for $\eta ^l_\lambda (x;q)$:
$$
\eta ^l_\lambda (x;q)=\frac{(yq^{l+1/2}e^\theta;q)_\infty
(yq^{l+1/2}e^{-\theta};q)_\infty}{({\rm i}yq^{1/2};q)_\infty
(-{\rm i}yq^{1/2};q)_\infty} =  \frac{(e^{{\rm i}\psi}q^{(2l+3)/4}x
\,e^\theta;q)_\infty (e^{{\rm i}\psi}q^{(2l+3)/4}xe^{-\theta};q)_\infty}
{({\rm i}e^{{\rm i}\psi}q^{(3-2l)/4}x;q)_\infty (-{\rm i}e^{{\rm i}\psi}
q^{(3-2l)/8}x;q)_\infty} .
$$

Employing here the $q$-binomial theorem
$$
\sum _{n=0}^\infty \frac{(a;q)_n}{(q;q)_n}\,t^n =
\frac{(at;q)_\infty}{(t;q)_\infty},\ \ \ \ |t|<1,\ \ |q|<1,
$$
and comparing the obtained expression with formula (8.5), we
deduce that
$$
P'_k(\lambda)=(-{\rm i})^kq^{k/2}\sum _{n=0}^k \frac{(-1)^n(-{\rm i}
\,q^le^\theta ;q)_n ({\rm i}q^le^{-\theta};q)_{k-n}}{(q;q)_n(q;q)_{k-n}} .
$$
It is a polynomial in $\cosh \theta =-\frac 12 dq^{-l+1/2}\lambda$. This
expression can be represented in terms of a basic hypergeometric functions.
To this end we take into account the relation
$$
(a;q)_{k-n}=\frac{(a;q)_kq^{n(n+1)/2}}{(-aq^k)^n(q^{1-k}a^{-1};q)_n}.
$$
As a result, we obtain
$$
P'_k(\lambda)=c'\sum _{n=0}^k\frac{(-{\rm i}q^le^\theta ;q)_n
(q^{-k};q)_n}{(q;q)_n (-{\rm i}q^{-k-l+1}e^\theta ;q)_n}
({\rm i}q^{-l+1}e^\theta )^n
$$  $$
= c'{}_2\phi_1(-{\rm i}q^le^\theta ,\
q^{-k};\ -{\rm i}q^{-k-l+1}e^\theta ;\ q,{\rm i}q^{-l+1}e^\theta ),
$$
where $c'=(-{\rm i}q^{1/2})^k({\rm i}q^le^{-\theta};q)_k/(q;q)_k$.
Applying here the relation (III.8) and then (I.8) from Ref. 24
we obtain the expresion
$$
P'_k(\lambda)=(-{\rm i})^kq^{-k(k+l-3/2)}\frac{(q^{2l};q)_k}
{(q;q)_k}\, {}_3\phi_1
(q^{-k},-{\rm i}q^le^\theta ,-{\rm i}q^le^{-\theta};\,
q^{2l};\; q,-q^k)  \eqno (8.7)
$$
for $P'_k(\lambda)$,  which explicitly exhibits a polynomial
dependence on $\cosh \theta$.  

These polynomials in $\lambda$ have many orthogonality relations, since a
closure of the operator ${I}_4^{(\psi)}$ is not selfadjoint. It is
difficult problem to find orthogonality relations explicitly. For this reason,
we cannot construct selfadjoint extensions of the operator ${I}_4^{(\psi)}$
and their spectra. Note that the polynomials (8.7) are very similar to
(but not coinciding with) the polynomials (5.17) in Ref. 36.

As in the previous case, we can realize the representation $T^+_l$
on the linear space of all polynomials $p(\lambda)$ in $\lambda$,
such that
$$
{I}_4^{(\psi)} p(\lambda)=\lambda p(\lambda).
$$

The polynomials $P'_k(\lambda)$, $\noindent
k=0,1,2,\cdots$ (and also the
polynomials $P_k(\lambda)$, $k=0,1,2,\cdots$), form a basis in
this space. It follows from (8.2) and (8.3) that
$$
{I}_4^{(\psi)}P_k(\lambda)= q^{-l-n-1/2}\left( e^{{\rm i}
\psi}\sqrt{(1-q^{k+1})(1-q^{2l+k})} \,P_{k+1}(\lambda)  \right.
$$   $$  \left.
+e^{-{\rm i}\psi}\sqrt{(1-q^{k})(1-q^{2l+k-1})}\,P_{k-1}(\lambda)\right) ,
$$
that is, ${I}_4^{(\psi)}$ acts on the basis $P_k(\lambda)$, $k=0,1,2,\cdots$,
by the same formula as it acts on the basis (2.15).
\bigskip

\noindent{\sf ACKNOWLEDGMENTS}

\medskip

This research has been supported in part by the SEP-CONACYT 
project 41051-F. A.U. Klimyk acknowledges the Consejo Nacional 
de Ciencia y Technolog\'{\i}a (M\'exico) for a C\'atedra 
Patrimonial Nivel II. We are grateful to J.S.Christiansen for 
the discusion of problems, related with extremal orthogonality 
measures for the $q$-Laguerre polynomials. 
\bigskip

\noindent{\sf REFERENCES}

\medskip

\noindent
${}^1$S. Lang, $SL_2(R)$ (Addison-Wesley, Massachusetts, 1975).

\noindent
${}^2$N. Ja. Vilenkin and A. U. Klimyk, {\it Representation of Lie
Groups and Special Functions}, vol. 1 (Kluwer, Dordrecht, 1991).

\noindent
${}^3$R. Howe and E. C. Tan, {\it Non-Abelian Harmonic Analysis.
Application of $SL(2,{\Bbb R})$} (Springer, New York, 1992).

\noindent
${}^4$D. Basu and K. B. Wolf, ``The unitary irreducible
representations of $SL(2,R)$ in all subgroup reductions'',
J. Math. Phys. {\bf 23}, 189--205 (1982).

\noindent
${}^5$J. Van der Jeugt, ``Coupling coefficients for Lie algebra
representations and addition formulas for special functions'',
J. Math. Phys. {\bf 38}, 2728--2740 (1997).

\noindent
${}^6$H. T. Koelink and J. Van der Jeugt, ``Convolutions for
orthogonal polynomials from Lie and quantum algebra
representations'', SIAM J. Math. Anal. {\bf 29}, 794--822 (1998).

\noindent
${}^7$J. Van der Jeugt and R. Jagannathan, ``Realizations of
$su(1,1)$ and $U_q(su(1,1))$ and generating functions for
orthogonal polynomials'', J. Math. Phys. {\bf 39},
5062--5078 (1998).

\noindent
${}^8$H. T. Koelink and J. Van der Jeugt, ``Bilinear generating
functions for orthogonal polynomials'', Constructive  Approximation
{\bf 14}, 481--497 (1999).

\noindent
${}^9$W. Groenevelt and E. Koelink, ``Meixner functions and
polynomials related to Lie algebra representations'', J. Phys. A:
Math. Gen. {\bf 35}, 65--85 (2002).

\noindent
${}^{10}$W. Groenevelt, E. Koelink, and H. Rosengren, ``Continuous
Hahn functions as Clebsch--Gordan coefficients'', math.CA/0302251.

\noindent
${}^{11}$V. G. Drinfeld, ``Hopf algebras and quantum Yang--Baxter
equation'', Soviet Math. Dokl. {\bf 32}, 254--258 (1985).

\noindent
${}^{12}$M. Jimbo, ``A $q$-analogue of $U( g)$ and the Yang--Baxter
equation'', Lett. Math. Phys. {\bf 10}, 63--69 (1985).

\noindent
${}^{13}$N. Yu. Reshetikhin, L. A. Takhtajan, and L. D. Faddeev, 
``Quantization of Lie groups and Lie algebras'', Leningrad Math. J.
{\bf 1}, 193--225 (1990).

\noindent
${}^{14}$A. U. Klimyk and I. I. Kachurik, ``Spectra, eigenvalues and
overlap functions for representation operators of $q$-deformed
algebras'', Commun. Math. Phys. {\bf 175}, 89--111 (1996).

\noindent
${}^{15}$H. Rosengren, ``A new quantum algebraic interpretation of the
Askey--Wilson polynomials'', Contemp. Math. {\bf 254},
371--394 (2000).

\noindent
${}^{16}$E. Koelink and J. V. Stokman, ``Fourier analysis on the
quantum $SU(1,1)$ group'', Publ. RIMS {\bf 37}, 621--715 (2001).

\noindent
${}^{17}$L. L. Vaksman and L. I. Korogodskii, ``Spherical functions
on the quantum group $SU(1,1)$ and a $q$-analogue of the
Mehler--Fock formula'', Funct. Anal. Appl. {\bf 25}, 48--49 (1991).

\noindent
${}^{18}$T. Kakehi, T. Masuda, and K. Ueno, ``Spectral analysis of a
$q$-difference operator which arises from the quantum $SU(1,1)$ group'',
J. Operator Theory {\bf 33}, 159--196 (1995).

\noindent
${}^{19}$A. Ballesteros and S. M. Chumakov, ``On the spectrum of a
Hamiltonian defined on ${\rm su}_q(2)$ and quantum optical
models'', J. Phys. A: Math. Gen. {\bf 32}, 6261--6269 (1999).

\noindent
${}^{20}$M. Noumi, T. Umeda, and M. Wakayama, ``Dual pairs, spherical
harmonics and a Capelli identity in quantum group theory'',
Compos. Math. {\bf 104}, 227--277 (1996).

\noindent
${}^{21}$N. Z. Iorgov and A. U. Klimyk, ``The $q$-Laplace operator
and $q$-harmonic polynomials on the quantum vector space'', J. Math.
Phys. {\bf 42}, 1326--1345 (2001).

\noindent
${}^{22}$N. Z. Iorgov and A. U. Klimyk, ``A Laplace operator
and harmonics on the quantum complex vector space'', J. Math.
Phys. {\bf 44}, 823--848 (2003).

\noindent
${}^{23}$N. M. Atakishiyev and A. U. Klimyk, ``Diagonalization of
representation operators for the quantum algebra $U_q({\rm
su}_{1,1})$'', Methods of Functional Analysis and Topology
{\bf 8}, No. 3, 1--12 (2002).

\noindent
${}^{24}$G. Gasper and M. Rahman, {\it Basic Hypergeometric Functions}
(Cambridge University Press, Cambridge, 1990).

\noindent
${}^{25}$G. E. Andrews, R. Askey, and R. Roy, {\it Special Functions}
(Cambridge  University Press, Cambridge, 1999).

\noindent
${}^{26}$I. M. Burban and A. U. Klimyk, ``Representations of the
quantum algebra $U_q({\rm su}_{1,1})$'', J. Phys. A: Math. Gen.
{\bf 26}, 2139--2151 (1993).

\noindent
${}^{27}$T. H. Koornwinder, ``Askey--Wilson polynomials as zonal spherical
functions on the $SU(2)$ quantum group'', SIAM J. Math. Anal.
{\bf 24}, 795--813 (1993).

\noindent
${}^{28}$N. M. Atakishiyev and P. Winternitz, ``Bases for
representations of quantum algebras'', J. Phys. A: Math. Gen.
{\bf 33}, 5303--5313 (2000).

\noindent
${}^{29}$N. M. Atakishiyev and A. U. Klimyk, ``Diagonalization of
operators and one-parameter families of nonstandard bases for
representations of ${\rm su}_q(2)$'', J. Phys. A: Math. Gen.
{\bf 35}, 5267--5278 (2002)..

\noindent
${}^{30}$R. Koekoek and R. F. Swarttouw, {\it The Askey-Scheme of
Hypergeometric Orthogonal Polynomials and Its $q$-Analoque}, Delft
University of Technology Report 98--17; available from {\tt ftp.tudelft.nl}.

\noindent
${}^{31}$Ju. M. Berezanskii, {\it Expansions in Eigenfunctions of
Selfadjoint Operators} (American Mathematical Society, Providence, RI, 1968).

\noindent
${}^{32}$T. S. Chihara and M. E. H. Ismail, ``Extremal measures for a
system of orthogonal polynomials'', Constructive Approximation {\bf 9},
111--119 (1993).

\noindent
${}^{33}$J. S. Christiansen, ``The moment problem associated with 
the $q$-Laguerre polynomials'', Constructive Approximation  {\bf 19},
1--22 (2003).

\noindent
${}^{34}$N. Cicconi, E. Koelink, and T. H. Koornwinder, ``$q$-Laguerre 
polynomials and big $q$-Bessel functions and their orthogonality
relations'',  Methods of Applied Analysis {\bf 6}, 109--127 (1999).

\noindent
${}^{35}$N. I. Akhiezer and I. M. Glazman, {\it The Theory of Linear
Operators in Hilbert Spaces} (Ungar, New York, 1961).

\noindent
${}^{36}$Ch. Berg and M. E. H. Ismail, ``$q$-Hermite polynomials and 
classical orthogonal polynomials'', Canad. J. Math. {\bf 48},
43--63 (1996).

\end{document}